\documentclass[12pt,aos]{article}
\usepackage{mathrsfs}
\usepackage{epsf}
\usepackage{amsthm,amsmath,amssymb,graphicx,graphics,epsfig}
\usepackage{lineno}
\usepackage{hyperref}
\usepackage{booktabs}
\usepackage[authoryear,round]{natbib}
\usepackage{longtable}
\usepackage{float}
\usepackage{multirow}
\usepackage{graphicx}
\usepackage[graphicx]{realboxes}
\usepackage{rotating}
\usepackage{algorithm}
\usepackage{color}
\usepackage{algorithmic}
\usepackage{amsmath}
\def\var{\mbox{var}}

\def\sign{\mbox{sign}}

\def\sgn{\mbox{sgn}}

\def\bSig\mathbf{\Sigma}

\newcommand{\Rmnum}[1]{\uppercase\expandafter{\romannumeral #1}}

\begin{document}
\baselineskip 8mm \setcounter{page}{0} \thispagestyle{empty}
\begin{center}
{\Large \bf High-dimensional outlier detection and variable selection via adaptive weighted mean regression } \vspace{3mm}

Jiaqi Li$^{1}$, Linglong Kong$^{1}$, Bei Jiang$^{1}$ and Wei Tu$^{2^*}$\\
{\it $^1$Department of Mathematical and Statistical Sciences, University of Alberta, Edmonton, AB, Canada \\
	$^{2^*}$ Department of Public Health Sciences and Canadian Cancer Trials Group, Queen's University, Kingston, ON, Canada}  \\
$^{*}${\it Corresponding author. Email address: wei.tu@queensu.ca}

\end{center}

\noindent{\bf Abstract} 
This paper proposes an adaptive penalized weighted mean regression for outlier detection of high-dimensional data.
In comparison to existing approaches based on the mean shift model, the proposed estimators demonstrate robustness against outliers present in both response variables and/or covariates.
By utilizing the adaptive Huber loss function, the proposed method is effective in high-dimensional linear models characterized by heavy-tailed and heteroscedastic error distributions.
The proposed framework enables simultaneous and collaborative estimation of regression parameters and outlier detection.
Under regularity conditions, outlier detection consistency and oracle inequalities of robust estimates in high-dimensional settings are established.
Additionally, theoretical robustness properties, such as the breakdown point and a smoothed limiting influence function, are ascertained. 
Extensive simulation studies and a breast cancer survival data are used to evaluate the numerical performance of the proposed method, 
demonstrating comparable or superior variable selection and outlier detection capabilities.
 

\newpage
\section{Introduction}
Outlier detection is a common task among statistical data analysis. Extensive research have been devoted to this topic
\cite{Rou03,maro06}, but outlier detection for high-dimensional data still remains challenging. 
\cite{zhao13} proposed a high dimensional influence
measure (HIM) to flag those points that greatly influence the calculated value of
marginal correlations as influential data. However, the idea of leave-one-out observation is not applicable for detecting
multiple outliers  with strong masking and/or swamping effects \cite{hadi93}. Methods based on projection pursuit have been proposed for
identifying high-dimensional multivariate outliers \cite{prieto01,maro02}, but they are usually computationally intensive. Recently, \cite{zhao19} proposed a novel group deletion 
procedure for multiple influential point detection using min- and max- extreme statistics based on a marginal-correlation-based influence measure, which effectively 
overcomes the masking and/or swamping effects.

Based on the mean-shift model, \cite{she11} considered outlier detection using non-convex penalized
regression on the mean-shift parameters, and the threshold-based iterative procedure was linked with M-estimators. \cite{kong18} 
investigated the problem of detecting outliers and selecting significant predictors in high-dimensional linear regression. Specifically, the nonzero mean-shift parameter estimates were identified as outliers and an $L_1$ penalty was imposed on the regression parameters to facilitate the selection of important predictors.
\cite{karu19} proposed a new penalized procedure that simultaneously attains full efficiency and maximum robustness by imposing adaptive weights on both the decision loss and penalty function.
Furthermore, the optimal convergence rate of $L_1$-penalized Huber's M-estimator under the mean-shift model
was studied by \cite{dala19}. 

Robustness can also be achieved by down-weighting those leverage points
if they can be detected in advance. The penalized weighted least absolute
deviation regression (PWLAD) method proposed by \cite{gao18} 
concurrently estimated the weights and regression coefficients of each observation via a lasso-type penalty on the weight vector.
Later, \cite{jiang20} extended it to achieve variable selection and outlier detection simultaneously.
Since the initial weight they used is based on MCD estimator, it is not suitable for sparse models. Additionally, the efficiency of estimators is compromised when dealing with heavy-tailed data. To address these limitations, we propose an adaptive penalized weighted mean regression approach that facilitates outlier detection and variable selection. Moreover, the proposed approach provides robustness against multiple outliers and heavy-tailed distributions. Notably, the proposed method accommodates both homogeneity and heterogeneity in the distributions of individual observations.


In classical statistical literature, various methods have been proposed to measure the robustness of estimation procedures.
Breakdown points and influence functions of estimators are commonly used to evaluate the robustness of existing methods.
In high dimension setting, \cite{wang13} analyzed the breakdown point of a certain nonconvex penalized M-estimator, while theoretical properties of regularized robust M-estimators have been studied by \cite{loh17}. Although influence functions of many penalized M-estimators
have been derived to evaluate the robustness of their estimators (see \cite{wangxq13}), and a rigorous framework has been established by \cite{marco17} using the 
limiting influence function for penalized M-estimators in the sparse model, one of the key assumptions of them is that the loss function has bounded second-order derivatives, which is not applicable to the Huber's loss function. To our best of our knowledge, the explicit form of influence function for Huber's estimator has not been established, and the robustness for adaptive Huber regression \cite{sun20} has not been evaluated using the influence function. To fill these gaps, we derive a smoothed limiting influence function for the adaptive Huber's estimator to evaluate its robustness. This technique can also be extended to any other loss functions without second derivatives.

The rest of the article is structured as follows: In Section 2, we construct a penalized weighted Huber-LASSO (PWHL) estimator for variable selection and outliers detection.
In Section 2.1, we implement an iterative algorithm to solve the bi-convex optimization problem, with criteria for tuning parameter selection being established in Section 2.2.
Theoretical properties of the proposed estimators are explored in Section 3, with an evaluation of robustness based on the breakdown point and the smoothed limiting influence function.
In Section 4, we carry out comprehensive simulation studies to assess the performance of the proposed methodology. 
An analysis of breast cancer data using the proposed technique is presented in Section 5. Lastly, in Section 6, we draw some conclusions.

\section{Methodology}
Let $\left(x_{i}, y_{i}\right)$ represent a random sample drawn from population $(X,Y)$ for the $i$th subject, where $i=1,\cdots, n$. 
Assume $Y$ is an univariate response variable, and $x_i$ is a $p_n$-dimensional predictor. Suppose $(x_i,y_i)$ satisfy a linear regression model
\begin{equation}\label{linear}
	y_{i}=x_{i}^{\top} \boldsymbol{\beta}+\varepsilon_{i},
\end{equation}
where $\boldsymbol{\beta}$ is a $p_n$-dimensional vector of unknown parameters, with the number of nonzero coefficients being $d_n$, namely $|\boldsymbol{\beta}|=d_n$. 
The error terms $\left\{\varepsilon_i, i=1,\ldots,n\right\}$ are i.i.d, which allows 
conditional heteroscedastic models that $\varepsilon_i$ can depend on $x_i$, concretely, $\varepsilon_{i}=\sigma\left(x_{i}\right) \tilde{\varepsilon}_{i}$ with 
function $\sigma(x_{i})$ of $x_i$, and $\tilde{\varepsilon}_{i}$ is independent of $x_i$. Assume the distributions of $X$ and $\varepsilon|X$ have zero mean, and 
consequently, $ \boldsymbol{\beta}$ is associated with the mean effect of $Y$, conditional on $X$.


To achieve variable selection and outlier detection in responses and/or covariates simultaneously, 
we consider the penalized weighted Huber-LASSO estimator:
\begin{equation}\label{obj}
	(\widehat{\boldsymbol{\beta}}, \widehat{\mathbf{w}})=\underset{\boldsymbol{\beta}, \mathbf{w}}{\arg \min }\left\{\frac{1}{2} \sum_{i=1}^{n} 
	\ell_{\alpha}\left(w_{i}(y_{i}-x_{i}^{\top} \boldsymbol{\beta})\right) +\mu \sum_{i=1}^{n} {\varpi}_{i}\left|1-w_{i}\right|+\lambda \sum_{j=1}^{p_n} \left|\beta_{j}\right|\right\},
\end{equation}
where the weights vector $\mathbf{w}=\left(w_{1}, \ldots, w_{n}\right)^{\prime}$, and $0<w_{i} \leq 1$ indicates the outlying effects for each observation, 
with $w_i<1$ representing a potential outlier and $w_i=1$ representing a normal observation.
Define the true values of $\boldsymbol{\beta}$ and $\mathbf{w}$ are $\boldsymbol{\beta}^{*}$ and $\mathbf{w}^{*}$ respectively.
Similar to \cite{gao18}, a prior weight $\varpi_{i}$ can be obtained by setting $\varpi_{i}=1 /\left|\log \left(w_{i}^{(0)}\right)\right|$ with an initial value
of weight $w^{(0)}_i$. When there are outliers in observations, $w_i$ can be used to downweight the influence of outliers.

The Huber loss \citep{huber64}
\begin{equation}\label{huber}
	\ell_{\alpha}(x)=\left\{
	\begin{array}{lll}
		2 \alpha^{-1}|x|-\alpha^{-2} & \text {if } & |x|>\alpha^{-1} \\
		x^{2} & \text { if } & |x| \leq \alpha^{-1}
	\end{array} \right.
\end{equation}
merges quadratic loss and least absolute deviation (LAD) by allowing the tuning parameter $\alpha$ to vary. As  
$\alpha$ converges to 0, it diminishes biases in estimating mean regression when the conditional distribution
of $\varepsilon_i$ is not symmetric \citep{fan17}.

\textbf{Computational algorithm:}
As objective function (\ref{obj}) exhibits bi-convexity in terms of $\mathbf{w}$ and $\boldsymbol{\beta}$, a solution path for $(\boldsymbol{\beta},\mathbf{w})$ 
can be derived by optimizing one variable while holding the other constant, akin to the algorithm proposed by\citet{jiang20}.
 A summary of the computational procedure for $\boldsymbol{\beta}$ and $\mathbf{w}$ is provided in the subsequent pseudocode.

\begin{algorithm}
	\caption{The solution path for Bi-optimization (2.2) given $\lambda$ and $\alpha$}
	\label{alg}
	\begin{algorithmic}
		\REQUIRE Initial estimates $\widehat{\boldsymbol{\beta}}^{(0)}$, $\widehat{\mathbf{w}}^{(0)}$, and thus $\varpi_{i}=1 /\left|\log \left(\widehat{w}_{i}^{(0)}\right)\right|$ for $i=1,\ldots,n$.
		\STATE Let $k=1$ and $\bar{\mu}_i=\mu \varpi_{i}$
		\WHILE {not converged}
		\STATE [Update $\boldsymbol{\beta}$] \\
		\STATE Denote $\ell^{\star}_{\alpha}\left(\boldsymbol{\beta}\right)=\ell_{\alpha}\left(w_{i}(y_{i}-x_{i}^{\top} \boldsymbol{\beta})\right)$
		\STATE Let 
		\begin{equation}\label{beta}
			\widehat{\boldsymbol{\beta}}^{(k)}=\arg \min _{\boldsymbol{\beta}}\left\{\sum_{i=1}^{n}\ell^{\star}_{\alpha}\left(\boldsymbol{\beta}\right)+\lambda \sum_{j=1}^{p_n} \left|\beta_{j}\right|\right\}
		\end{equation} 
		\STATE [Update $\mathbf{w}$]
		\STATE $r_i^{(k)}=y_i-x_i^{\top}\widehat{\boldsymbol{\beta}}^{(k)}$, for $i=1,\ldots,n$
		\STATE If $\ell_{\alpha}\left(r_{i}^{(k)}\right)>\bar{\mu}_{i}$, let $\mathbf{w}_{i}^{(k)} \leftarrow \bar{\mu}_{i} /\ell_{\alpha}\left(r_{i}^{(k)}\right)$, otherwise $\mathbf{w}_{i}^{(k)} \leftarrow 1$ 
		\STATE converged $\leftarrow\left\|\mathbf{w}_{i}^{(k)}-\mathbf{w}_{i}^{(k-1)}\right\|_{\infty}<\epsilon=10^{-8}$
		\STATE $k \leftarrow k+1$
		\ENDWHILE
		\STATE Output $\widehat{\boldsymbol{\beta}}=\widehat{\boldsymbol{\beta}}^{(k)}$ and $\widehat{\mathbf{w}}=\mathbf{w}^{(k)}$
	\end{algorithmic}
\end{algorithm}

The convex optimization problem (\ref{beta}) bears resemblance to RA-Lasso \citep{fan17}, and can be solved using gradient-descent-based algorithms.
By employing a Taylor expansion, the original loss function is replaced with a local isotropic quadratic approximation (LQA), which allows for the iterative minimization of an objective function.

 \text{\bf Remark 1.} While the optimization problem in Eq.(2.2) is fundamentally non-convex, for a given weight $\mathbf{w}$, the iterative updating procedure (2.4) can be viewed as a sequential univariate optimization w.r.t unknown parameter $\boldsymbol{\beta}$.
	This procedure represents a convex regularized optimization problem, as the weighted Huber loss function $\ell_\alpha\left(w_i\left(y_i-x_i^{\top} \boldsymbol{\beta}\right)\right)$ satisfies the local Restricted Strong Convexity (RSC) condition with high probability\citep{sun20}.
	When provided with the updated estimator $\widehat{\boldsymbol{\beta}}$, the $L_1$ regularized optimization problem concerning $\mathbf{w}$ also constitutes a convex optimization due to the convexity of the Huber loss function under a moderate sample size $n$. 
	The Alternating Convex Search\citep{Gor07} algorithm can be employed to address this biconvex optimization. A Greedy Coordinate Descent algorithm\citep{wu08} can be designed for the outer layer, updating the weight vector $\mathbf{w}$, while being nested with a gradient descent algorithm employing local quadratic approximation (LQA) for the inner layer, updating $\boldsymbol{\beta}$ in the context of solving (2.4).

 The sequence of joint optimal estimates $\left(\hat{\boldsymbol{\beta}}^{(k)}, \hat{\mathrm{w}}^{(k)}\right)$ for k-th iteration 
	possesses at least one limit point\citep{Gor07}. 
	Although it can be readily demonstrated that the limit point of such a sequence is a partial optimum, there is no guarantee that the convergence of the sequence must be a local optimum; however, it remains sufficiently close for practical purposes. Nonetheless, certain regularity conditions ensure that the iterative algorithm converges to a unique local minimizer. A comprehensive analysis is presented in Appendix B of the Supplementary Information.

\textbf{Selection of tuning parameters:}
In order to execute the aforementioned algorithm, the tuning parameter $\alpha$ and
a pair of regularization parameters $(\mu,\lambda)$ must be appropriately selected from a grid of candidates. This selection process influences outlier detection consistency and variable selection efficiency.

Given the significant impact of the tuning parameter $\mu$ on the performance of both outlier identification and parameter estimation, we employ a random weighting method to select an optimal $\widehat \mu$ from a 
grid of $\mu$s as recommended by \citet{gao18}. Denote random weights vector $\boldsymbol{\omega}=\left(\omega_{1}, \cdots, \omega_{n}\right)^{\prime}$, where weights 
$\omega_{1}, \cdots, \omega_{n}$ are i.i.d satisfying $E\left(\omega_{i}\right)=\operatorname{Var}\left(\omega_{i}\right)=1$.
By solving a perturbed objective function, we acquire the corresponding estimates in relation to a set of $\mu$ and $\boldsymbol{\omega}$ values from:
$$
\begin{aligned}
	\left(\widehat{\boldsymbol{\beta}}(\mu ; \boldsymbol{\omega}), \widehat{\mathbf{w}}(\mu ; \boldsymbol{\omega})\right)
	=  \underset{\boldsymbol{\beta},\mathbf{w}}{\arg \min }\left\{\frac{1}{2} \sum_{i=1}^{n}  \ell_{\alpha}\left(\omega_{i}{w_{i}}(y_{i}-x_{i}^{\top} \boldsymbol{\beta})\right)+\mu \sum_{i=1}^{n} {\varpi}_{i}\left|1-w_{i}\right|\right\}.
\end{aligned}
$$
For any two sets of random weights $\boldsymbol{\omega}_{1}$ and $\boldsymbol{\omega}_{2}$, there are two perturbed weight estimates: 
$\widehat{\mathbf{w}}\left(\mu ; \boldsymbol{\omega}_{1}\right)$ and $\widehat{\mathbf{w}}\left(\mu ; \boldsymbol{\omega}_{2}\right)$,
which lead to two sets of suspected outliers, denoted as $\mathcal{O}\left(\mu ; \boldsymbol{\omega}_{1}\right)$ and $\mathcal{O}\left(\mu ; \boldsymbol{\omega}_{2}\right)$.
Following \citet{gao18}, we utilize Cohen's kappa coefficient $\kappa\left(\mathcal{O}\left(\mu ; \boldsymbol{\omega}_{1}\right), \mathcal{O}\left(\mu ; \boldsymbol{\omega}_{2}\right)\right)$ to
represent the agreement between these two sets of suspected outliers, and we iteratively generate $B$ pairs of random weights $(\boldsymbol{\omega}_{b1},\boldsymbol{\omega}_{b2})$ for $b=1,\ldots, B$.
Denote
\begin{equation}\label{cohen}
	{S}(\mu)=\frac{1}{B} \sum_{b=1}^{B} \kappa\left(\mathcal{O}\left(\mu ; \boldsymbol{\omega}_{b 1}\right), \mathcal{O}\left(\mu ; \boldsymbol{\omega}_{b 2}\right)\right),
\end{equation}
Consequently, we choose the optimal $\widehat\mu$ corresponding to the minimizer of ${S}(\mu)$ across a spectrum of candidates. In our numerical studies, we seek the optimal $\widehat\mu$ through a grid of candidates, with values ranging from $0.1$ to $0.5$ with an interval of $0.1$.

Next, we apply the BIC-type criterion to select the tuning parameters $\alpha$ and $\lambda$:
\begin{equation}\label{bic}
	\operatorname{BIC}(\alpha,\lambda)=n \log (\operatorname{RSS} / n)+df \left\{ \log(n)+c*\log(p_n+n)\right\},
\end{equation}
where $df$ represents the number of non-zero components of $\widehat{\boldsymbol{\beta}}_{\alpha,\lambda}$ and $1-\widehat{\boldsymbol{w}}$,
and the residual $\mathrm{RSS}=\sum_{\mathrm{i}=1}^{\mathrm{n}} \hat{e}_{\mathrm{i}}^{2}$
with $\hat{e}_{i}^{2}=\left\{\widehat w_{i}\left(y_{i}-x_{i}^{\top} \widehat{\boldsymbol{\beta}}_{\alpha,\lambda}\right)\right\}^{2}$ with respect to
estimators $ \widehat{\boldsymbol{\beta}}_{\alpha,\lambda}$ and $\widehat{\mathbf{w}}$ of (\ref{obj}) for given $\alpha$, $\lambda$ and $\widehat \mu$.
 The constant prior $c\in(0,2)$ is dependent with $\kappa$(defined later) to ensure the selection consistency\citep{chen08},  
	and in our simulations, we set $c=1.01$.

We obtain the optimal ($\hat\alpha$, $\hat\lambda$) minimizing (\ref{bic}). 
When error terms are homogeneous, the candidate range for $\alpha$ is set within the interval $(0.1, 1.0)$ with increments of $0.1$.
In cases of heteroscedasticity, a smaller value range $(0.01, 0.8)$ is adopted.
We search for an optimal $\widehat\lambda$ between $0.1$ and $1.0$ at intervals of $0.1$. For outlier detection, we suggest $\alpha=0.1$, while for robust parameter estimation, $\alpha=0.01$ is recommended. This choice is computationally efficient in both numerical studies and real data analysis.

\section{Theoretical results}

\subsection{Outlier detection consistency}
Our theoretical analysis builds upon and enhances the work of \citet{gao18}, in which the outlier detection properties of the penalized weighted least absolute deviation were established. We relax constraints on noise and permit sub-Gaussian tails in the mean regression, thus extending the applicability of the model.

We follow similar notations as \citet{gao18}, where the true outlier set is
$\mathcal{O}=\left\{1 \leq i \leq n: 0<w_{i}^{*}<1\right\}$ with $|\mathcal{O}|=q_n$, representing the number of outlier points.
The estimated set is $\widehat{\mathcal{O}}=\left\{1 \leq i \leq n: 0<\widehat{w}_{i}<1\right\}$. The weight satisfies $0<w_{i}^{*}\leq 1$, and
identifies outliers as those with values smaller than 1.
Denote $l=\max _{1 \leq i \leq n}\left\|x_{i j}\right\|_{\infty}$. 
In order to demonstrate that our proposed method can detect all true outliers with high probability, irrespective of the covariates' dimension $p_n$, we establish the following conditions:

C1. Assume $\varepsilon_i$ for $i=1,\ldots,n$ are i.i.d with sub-Gaussian distribution with zero mean and
variance proxy $\sigma^2$, namely, $\varepsilon \sim \operatorname{subG}\left(\sigma^{2}\right)$. 

C2. There exists $\underline{\varpi}_{n}>0$ and $\bar{\varpi}_{n}>0$ satisfying 
$$
P\left(\left\{\max _{i \in \mathcal{O}} \varpi_{i} \leq \underline{\varpi}_{n}\right\} \cap\left\{\min _{i \in \mathcal{O}^{c}} \varpi_{i} \geq \bar{\varpi}_{n}\right\}\right)=1-o(1),
$$

C3. For tuning parameters $\alpha$ and $\mu$ together with $q_n$, $\underline{\varpi}_{n}$, and $\bar{\varpi}_{n}$, they satisfy:\\
(i) $\underline{\varpi}_{n}=o\left(\left(q_{n} \mu \alpha \right)^{-1}\right)$, $\alpha^{-1}q_n=o(1)$,
(ii) $1 / \bar{\varpi}_{n}=o\left(\alpha\mu / \sqrt{\log (n)}\right)$.\\ 

\text{\bf Remark 2.} In Condition C1, the error distribution of $\varepsilon_i$ is not limited to a continuous and positive density
at the origin, which is a requirement in \cite{gao18}.  
More broadly, we permit 
$\varepsilon_i$ to be conditionally i.i.d with $\operatorname{subG}\left(\sigma^{2}\right)$ given $x_i$ in heterogeneous model, and same conclusion can be drawn since 
the conditional tails of $\varepsilon_i$ can be controlled under regularity conditions. Conditions C2 and C3 constraint $\varpi_{i}$ to be small in the true outlier set, while being large in normal data, as suggested in\citet{gao18}. In contrast, C3 implies a higher order rate on $\underline{\varpi}_{n} / \bar{\varpi}_{n}$ than that in \citet{gao18}. Furthermore, 
C3 also implies that $\alpha$ must be sufficiently large as $q_n$ grows with $n$, which follows a higher order than $O(q_n)$, but lower than $O(\sqrt{\log n})$. 

{\bf Theorem 1.} Suppose Conditions C1--C3 hold, and a given initial estimate $\widehat{\boldsymbol{\beta}}^{(0)}$ satisfies
$$
P\left(\left\|\widehat{\boldsymbol{\beta}}^{(0)}-\boldsymbol{\beta}^{*}\right\|_{2}>\frac{\alpha}{2}\mu\max \left\{\bar{\varpi}_{n}/2, \underline{\varpi}_{n}\right\} /\left(\sqrt{p_n} l \right)\right) =o(1),
$$
for the estimate $\widehat{\mathbf{w}}$ of (\ref{obj}), 
we have 
$$
\lim _{n \rightarrow \infty} P(\widehat{\mathcal{O}}=\mathcal{O})=1.
$$ 

Theorem 1 suggests that the PWHL estimator can detect all outliers with a probability tending to 1 asymptotically when initial estimator $\widehat{\boldsymbol{\beta}}^{(0)}$ satisfies 
oracle property in ultra-high dimension. Such a consistent estimator can be obtained from penalized quantile regression \citep{wang11}.
In numerical studies, we utilize the sparse least trimmed squares estimator \citep{Alf13} as a warm start, 
which can be obtained using R function ``sparseLTS". For the initial estimate $\widehat{\mathbf{w}}^{(0)}$ in the 
weight $\varpi_{i}=1 /\left|\log \left(\widehat{w}_{i}^{(0)}\right)\right|$,
\citet{gao18} and \citet{jiang20} computed leverage values and assessed
data contamination severity based on the robust Mahalanobis distance, as proposed in \citet{fil08}. However, this approach is effective only when $p_n \leq n$, and 
lacks efficiency when strong masking/swamping effects are present.
To compute more efficient initial weights, we employ the MIP proposed by \citet{zhao19} to obtain initial values, which overcomes both masking and swamping effects.

\subsection{Non-asymptotic oracle inequalities}

In this section, we investigate the PWHuber-LASSO estimator in high dimensions, permitting $p_n$ to grow exponentially with the sample size $n$, 
specifically $\log p_n=O(n^{\kappa})$ for some $\kappa>0$.
We provide $L_1$-norm and $L_2$-norm oracle inequalities under certain mild conditions. We first introduce some notations.

Let $r_{i, \boldsymbol{\beta}}=\left(y_{i}-x_{i}^{\top} \boldsymbol{\beta}\right)/\varpi_{i}$, $v_i=\varpi_{i}(1-w_i)$. 
Denote the joint estimator $\boldsymbol{\theta}=\left(\boldsymbol{\beta}^{\top}, \frac{\mu}{\lambda}\boldsymbol{\nu}^{\top}\right)^{\top}$, where
$\boldsymbol{\nu}=(v_1,\ldots,v_n)^{\top}$. Thus we have $\lambda\|\boldsymbol{\theta}\|_{1}=\lambda\|\boldsymbol{\beta}\|_{1}+{\mu}\|\boldsymbol{\nu}\|_{1}$.
Analogously to $\boldsymbol{\beta}^{*}$, $\boldsymbol{\nu}^{*}$ and $\boldsymbol{\theta}^{*}$ are defined as the true vale of $\boldsymbol{\nu}$ and $\boldsymbol{\theta}$, respectively.
Let set $\mathcal{S}$ denote the active set of $\boldsymbol{\theta}^{*}$ and $\mathcal{S}_1=\{0\leq i\leq p_n:\beta_i^{*} \neq 0 \}$ with $|\mathcal{S}_1|=d_n$. Thus $\mathcal{S}=(\mathcal{S}_1,\mathcal{O})$,
and $s_n=|\mathcal{S}|=d_n+q_n$ denote the cardinality of the sparse $\boldsymbol{\theta}$.

Define  $z_{i, \boldsymbol{\beta}}^{\top}=\left(x_{i}^{\top}, 0, \ldots,\left(\lambda / {\mu}\right) r_{i, \boldsymbol{\beta}},\ldots, 0 \right)$ as a $(p_n+n)$-dimensional vector with
$(p_n+i)$-th element being $\left(\lambda / {\mu}\right) r_{i, \boldsymbol{\beta}}$, for $i=1,\ldots,n$,
and let $\mathbf{Z}_{\beta}=\left(\begin{array}{c}z_{1, \beta}^{\top} \\ \cdots \\ z_{n, \beta}^{\top}\end{array}\right)$. It is easy to find that
(\ref{obj}) is equivalent to
\begin{equation}\label{objnew}
	\widehat{\boldsymbol{\theta}}_{{\mu},\lambda}=\underset{\boldsymbol{\beta}, \mathbf{w}}{\arg \min } L_{\alpha}(\boldsymbol{\theta})+\lambda\|\boldsymbol{\theta}\|_{1},
\end{equation}
where $L_{\alpha}(\boldsymbol{\theta})=\ell_{\alpha}(Y-\mathbf{Z}_{\boldsymbol{\beta}} \boldsymbol{\theta})$.

We start with the following regularity conditions:

C4. Assume $x_i\in \mathbf{R}^{p_n}$ are i.i.d from a sub-Gaussian random vector $x$ for $i=1,\ldots,n$.
The regression errors $\varepsilon_i$ are also i.i.d from a sub-Gaussian random vector $\varepsilon$, and satisfy $\mathbb{E}\left(\varepsilon_{i} \mid x_{i}\right)=0$ and
$v_{i, \delta}=\mathbb{E}\left(\left|\varepsilon_{i}\right|^{1+\delta} \mid x_{i}\right)<\infty$ for some $\delta>0$.
Without loss of generality, we assume $\varepsilon \sim \operatorname{subG}\left(\sigma^{2}\right)$.
Let $v_{\delta}=\frac{1}{n} \sum_{i=1}^{n} v_{i, \delta}<\infty$ for some $0<\delta \leq 1$, and $\nu_{\delta}=\min \left\{v_{\delta}^{1 /(1+\delta)}, v_{1}^{1 / 2}\right\}$.

C5. Let $\Sigma_n=\frac{1}{n}\sum_{i=1}^nz_{i,\beta^{*}}z_{i,\beta^{*}}^{\top}$, and $\Sigma=E(\Sigma_n)=(\sigma_{ij})_{1\leq i,j\leq (p_n+n)}$ is positive definite. Denote $\kappa_{l}$
as the minimum eigenvalues of $\Sigma$.

C6. Assume there exists a positive constant $M>0$ satisfies 
$\mathbb{P}(|\langle\boldsymbol{u}, \widetilde{z}_{i,\beta^{*}} \rangle| \geq t) \leq 2 \exp \left(-t^{2}\|\boldsymbol{u}\|_{2}^{2} / M^{2}\right)$ for all
$t \in \mathbb{R}$ and $\boldsymbol{u} \in \mathbb{R}^{p_n+n}$, where $\widetilde{z}_{i,\beta^{*}}=\boldsymbol{\Sigma}^{-1 / 2}z_{i,{\beta}^{*}} $. 

C7. For $\alpha_{0}^{-1} \geq \nu_{\delta}$. Let $\alpha^{-1}=\alpha_{0}^{-1}(n/\log(p_n+n))^{\max\{1/(1+\delta),1/2\}}$. 

C8. Let $\lambda/{\mu}\leq O(1)$, and $\lambda =2 L \sigma_{\max } \alpha_{0}^{-1}\{(\log(p_n+n)) / n\}^{\min \{\delta /(1+\delta), 1 / 2\}}$ with
$\sigma_{\max }=\max _{1 \leq j \leq p_n} \sigma_{j j}^{1 / 2}$ and $L=(2\sqrt{2}+1)M+1$.

\text{\bf Remark 3.} Conditions C4 and C5 are common in the literature in high-dimensional random design analysis, and they are similar to those in \citet{sun20}.
Condition C6 is easy to be satisfied since random vectors $z_{i, \boldsymbol{\beta}^{*}}$ for $i=1,\ldots,n$ are naturally i.i.d from a sub-gaussian distribution 
when Conditions 4 and 5 are both satisfied. 
For $\delta \geq 1$, Condition 7 indicates that the convergence rate of the adaptive robustification parameter $\alpha^{-1}$  
is at most in an order of $O(\sqrt{n/\log(p_n+n)})$. If $0<\kappa\leq 1$, there exists a contradiction with C3, where $\alpha^{-1}$ is advised to be small enough, creating a trade-off in consistency between outlier detection and parameter estimation. Conversely, $\kappa \geq 1$ signifies a higher order of $p_n$, and the consistency in outlier detection aligns with that in parameter estimation, provided a suitable $\alpha$ is chosen.

{\bf Theorem 2.} Assume Conditions C4--C8 hold, when $n \geq C\left(M, \sigma_{\max}, \kappa_{l}\right) s_n \log(p_n+n)$ where $C\left(M,\sigma_{\max}, \kappa_{l}\right)$ only depend on
$(M,\sigma_{\max},\kappa_{l})$. Then we give the upper bound for the $L_2$-error:
\begin{equation}\label{l2-error}
	\scriptsize{P\left(\left\|\Sigma^{1/2}\left(\widehat{\boldsymbol{\theta}}_{{\mu}, \lambda}-\boldsymbol{\theta}^{*}\right)\right\|_{2} \leq C \kappa_{l}^{-1 / 2} M \sigma_{\max } \alpha_{0}^{-1} (d_n+q_n)^{1 / 2}\left(\frac{\log (p_n+n)}{n}\right)^{\min \{\delta /(1+\delta), 1 / 2\}}\right)\geq 1-3 (p_n+n)^{-1}},
\end{equation}
for some constant $C\geq 0$.

We introduce the notation `$a_{n} \lesssim b_{n}$' for two sequences $\left\{a_{n}\right\}_{n \geq 1}$ and $\left\{b_{n}\right\}_{n \geq 1}$, which represents $a_{n} \leq K b_{n}$ for some constant $K>0$ independent of $n$.
Additionally, we present $L_1$-norm and $L_2$-norm oracle inequalities:
\begin{equation}\label{l1-norm}
	\small{\left\|\hat{\boldsymbol{\beta}}_{\mathcal{S}_{1}}-\boldsymbol{\beta}_{\mathcal{S}_{1}}^{*}\right\|_{1}+\left\|\hat{\mathbf{w}}_{\mathcal{O}}-\mathbf{w}_{\mathcal{O}}^{*}\right\|_{1} \lesssim \kappa_{l}^{-1} M \sigma_{\max } \alpha_{0}^{-1} (d_n+q_n)\left(\frac{\log(p_n+n)}{n}\right)^{\min \{\delta /(1+\delta), 1 / 2\}}},
\end{equation}
and 
\begin{equation}\label{l2-norm}
	\small{\left\|\hat{\boldsymbol{\beta}}_{\mathcal{S}_{1}}-\boldsymbol{\beta}_{\mathcal{S}_{1}}^{*}\right\|_{2}+\left\|\hat{\mathbf{w}}_{\mathcal{O}}-\mathbf{w}_{\mathcal{O}}^{*}\right\|_{2} \lesssim 
		\kappa_{l}^{-1 / 2} M \alpha_{0}^{-1} (d_n+q_n)^{1 / 2}\left(\frac{\log (p_n+n)}{n}\right)^{\min \{\delta /(1+\delta), 1 / 2\}}},
\end{equation}
with probability at least $1-3(p_n+n)^{-1}$.

Theorem 2 demonstrates that the $L_1$-norm errors and $L_2$-norm errors of PWHL estimator scales as $(d_n+q_n)\left(\frac{\log(p_n+n)}{n}\right)^{1/2}$ and $(d_n+q_n)^{1 / 2}\left(\frac{\log(p_n+n)}{n}\right)^{1/2}$ 
respectively, assuming that $\delta=1$ under regularity conditions.
Consistent estimators can be obtained with appropriately tuned robustification parameter $\alpha$ and regularization parameter $\lambda$ as suggested in \citet{sun20}.

\subsection{Robustness properties}

In this section, we evaluate the robustness properties of our estimator from a theoretical perspective. 

The influence function introduced by \citet{Ham71} measures the stability of estimators given an infinitesimal contamination.  
Suppose observations $\mathbf{Z}=(z_{1},\ldots,z_{n})$ are drawn from a common distribution $F$ over the space $\mathcal{Z}$, and the loss function $L : \mathcal{Z} \times \mathbb{R}^{p_n+n} \mapsto \mathbb{R}$
links the parameter space $\Theta \in \mathbb{R}^{p_n+n}$ with the observed data. The empirical distribution $\hat{F}=\frac{1}{n} \sum_{i=1}^{n} \delta_{z_{i}}$ 
is defined with the distribution probability
$\delta_{z}$ assigning mass 1 at the point $z$ and 0 elsewhere. Then the value $E_{\hat{F}}[L(Z, \theta)]$ can be an estimator 
of the unknown population risk function $E_{F}[L(Z, \theta)]$. 

To address the issue of nondifferentiable penalty functions,
the limiting influence function was proposed in \citet{marco17}, where a sequence of differentiable penalized M-estimators converging to
the penalized M-estimator of interest s used to derive the $\Lambda_{\lambda}$ satisfying two derivatives with respect to $\theta$. However, the loss function in \citet{marco17} is required to have second bounded derivatives, which is not suitable for the Huber loss function. To achieve more versatility, we propose a more general limiting influence function, where we relax the assumptions on derivative functions established in \citet{marco17}, and we allow both the loss function and penalty function to not necessarily have second-order derivatives.

Denote the first derivative of $\Lambda_{\lambda}$ as $U\left(Z;\boldsymbol{\theta}\right)$. For penalized adaptive Huber regression (\ref{objnew}), we have
\begin{equation}\label{ef}
	\begin{aligned}
		U\left(Z;\boldsymbol{\theta}\right)&=\psi_{\alpha}\left(Z, \boldsymbol{\theta}\right)+p^{\prime}(\boldsymbol{\theta} ; \lambda)\\
		&=\left(\alpha^{-1}\sign(e)I\{|e|>\alpha^{-1}\}+eI\{|e|\leq \alpha^{-1}\}\right)^{\top}\Delta_{\boldsymbol{\theta}}+p^{\prime}(\boldsymbol{\theta} ; \lambda)\\
		&\triangleq U_F\left(Z;\boldsymbol{\theta}\right)+p^{\prime}(\boldsymbol{\theta} ; \lambda),
	\end{aligned}
\end{equation}
where $e=Y-\mathbf{Z}_{\boldsymbol{\beta}} \boldsymbol{\theta}$,
$\psi_{\alpha}(\cdot)$ is the first derivative of (\ref{huber}), and $p^{\prime}(\cdot)$ is the first derivative of penalty function $p(\boldsymbol{\theta} ; \lambda)$ with respect to $\boldsymbol{\theta}$.
Here $\Delta_{\boldsymbol{\theta}}=\left[\rm -vec{(\mathbf{w}\cdot X)},\rm blockdiag\left\{-\frac{\lambda}{\mu \varpi_{i}}\left(y_i-x_i^{\top}\boldsymbol{\beta}\right)\right\}_{i=1,\ldots,n}\right]$, where
$\mathbf{w}\cdot X=(w_1x_1^{\top},\ldots,w_nx_n^{\top})$. 

Since the Huber loss function (\ref{huber}) only has the first derivative function, which contains non-differentiable indicator functions such as $I\{|e|>\alpha^{-1}\}$.
Motivated by the smoothed function proposed by \citet{heller07}, we define a smoothed function of $T(F)$ as $\widetilde T(F)$ in the form of:
\begin{equation}\label{sTF}
	\begin{aligned}
		\widetilde T({F})&=\operatorname{argmin}_{\theta} \widetilde{\Lambda}_{\lambda}(\theta ; F)\\
		&=\operatorname{argmin}_{\theta} \widetilde{L}_{F}\left(Z;\theta\right)+p(\theta ; \lambda),
	\end{aligned}
\end{equation} 
where $\widetilde{L}_{F}\left(Z;\theta\right)$ is a smoothed loss function, obtained by plugging in the differentiable replacement of indicator functions using an appropriate local distribution function $\Phi(\cdot)$, concretely, $I(z>\alpha^{-1})=\Phi\left(\frac{z-\alpha^{-1}}{h}\right)$ with a bandwidth $h$. 
Denote $\widetilde{\psi}_{\alpha}\left(Z, \boldsymbol{\theta}\right)$ as the smoothed version of $\psi_{\alpha}\left(Z, \boldsymbol{\theta}\right)$.
Let $\widetilde{\mathrm{IF}}(Z ;  F,  \widetilde T)$ be the influence function of $\widetilde T({F})$.

Under some mild conditions, the following theorem helps establish the limiting influence function of (\ref{objnew}), and we can demonstrate that it is asymptotically equivalent 
to the original influence function.

The proof of the theorem requires the following conditions:

C9. Assume the parameter vector $\boldsymbol{\theta}_{\alpha}^{*}$ belongs to an $L_q$-ball with a uniform radius $R_q$, namely 
$\sum_{j=1}^{p_n+n}\left|\theta_{\alpha,j}^{*}\right|^{q} \leq R_{q}$ for some $q\in (0,1]$, with $\boldsymbol{\theta}^{*}_{\alpha}$ converging to $\boldsymbol{\theta}^{*}$ as 
$\alpha$ goes to 0. The covariate vector $x_i$ satisfies $L=\max_{i}\|x_i\|_{\infty}<\infty$, and $R=\sup_{i}\|r_{i,\boldsymbol{\beta}}\|_{\infty}<\infty$.
Let $\lambda/{\mu}\leq O(1)$ in (\ref{objnew}), namely, $\left|\frac{\lambda}{{\mu}}\right|\leq C$ for some $C>0$.

C10. The local distribution function $\Phi(z)$ is continuous and its first derivative $\phi(z)=\partial \Phi(z) / \partial z$ is symmetric about zero with $\int z^{2} \phi(z)<\infty$, e.g.
the distribution function of the standard normal distribution.

C11. The bandwidth $h$ is chosen such that $n h \rightarrow 0$ and $n h^3 \rightarrow 0$, and the tuning parameter $\alpha^{-1}=\frac{1}{2}h\log(n)\rightarrow 0$, as $n \rightarrow \infty$.

\text{\bf Remark 4.} Condition C9  imposes constraints on residuals $r_{i,\boldsymbol{\beta}}$ within a bounded rectangular domain and bounds random covariate vectors $x_{i}$. This leads to
a uniformly bounded random matrix $\Delta_{\boldsymbol{\theta}}$. To be more specific, we 
have $\|\Delta_{\boldsymbol{\theta}}\|_{1}\leq \|\rm vec(X), C_1 R \rm I_{n \times n}\|_{1} \triangleq\|\bar{\bar{\Delta}}_{\boldsymbol{\theta}}\|_{1}$ for some positive real value $C_1>0$.

Condition C11 provides the convergence rate of $h$ with respect to the sample size $n$, and $\alpha^{-1}$ is required to converge to 0 as $n$ increases. In practice, our analysis of local robustness using the influence function does not depend on the bandwidth $h$, instead, the robustification parameter $\alpha$ dominates with a properly chosen $h$ that satisfies regularity conditions.     

{\bf Theorem 3.} Assume C9--C11 hold, consider a sequence $\left\{p_{m}\right\}_{m \geq 1} \in C^{\infty}(\Theta)$ converge to $p(\boldsymbol{\theta} ; \lambda) $ when $m \rightarrow \infty$.
For the approximating problems $\widetilde{\Lambda}_{\lambda}\left(\theta ; F, p_{m}\right)=E_{F}[\widetilde{L}(Z, \theta)]+p_{m}(\theta ; \lambda)$, we have:

(1) $\lim _{m} \widetilde T\left({F} ; p_{m}\right)=\widetilde T({F})$,

(2) the limiting influence function defined as $\widetilde{\mathrm{IF}}(z ; {F}, \widetilde T):=\lim _{m \rightarrow \infty} \widetilde{\mathrm{IF}}_{p_{m}}(z ; {F}, \widetilde T)$ does not depend on
the choice of $p_m$,

(3) the smoothed influence function $\widetilde{\mathrm{IF}}(z ; F, \widetilde{T}) \xrightarrow[n \rightarrow \infty]{P} \mathrm{IF}(z ; F, T)$,

(4) the influence function $\widetilde{\mathrm{IF}}(z ; F, \widetilde{T})$ has the following form:
$$
\widetilde{\mathrm{IF}}(z ; F, \widetilde T)=-\Sigma^{-1}\left(\widetilde{U}_F\left(z, \theta\right)+p_{\lambda}^{\prime}\left(\left|\theta\right|\right)\right),
$$
where $\Sigma^{-1}= \rm blockdiag \left\{\left(M_{11}+p^{\prime \prime}_{\lambda}\right)^{-1}, 0\right\}$, with $M_{11}=E_{F}\left[\dot{\widetilde{U}}_{F_{11}}(Z, \widetilde T(F))\right]$ and
$\rm p^{\prime \prime}_{\lambda}$ is a diagonal matrix with diagonal elements $p_{\lambda, j}^{\prime \prime}\left(\left|\theta_{j}\right|\right)$ for $j=1,\ldots,s_n$.

Theorem 3 demonstrates that the smoothed limiting functional $\widetilde T({F})$ is unique, and the limiting influence function 
does not depend on the choice of smooth penalty functions $p_{m}$. We derive the influence function and finite-sample breakdown point in more detail in the Appendix C. 
 Since the estimator in (\ref{obj}) is equivalent to a canonical penalized estimator with Huber loss after
	transformation, the PWHL estimator has a breakdown point of $\frac{1}{n}$, which has been implied in \citet{Alf13}. 
	Though we conclude such an undesirable result, our proposed estimator exhibits superior robustness against outliers and heteroskedasticity 
	from extensive simulation results compared to the estimator proposed in \citet{kong18}, which enjoys a high breakdown point that can be as high as $\frac{1}{2}$.


\section{Numerical studies}
We conduct simulation studies to assess the performance of our proposed method PWHuber-LASSO on the linear model (\ref{linear}), and compare it with the latest procedures: PM proposed by
\citet{kong18}, AHuber proposed by \citet{sun20}, which can be computed using the function ``tfNcvxHuberReg" in R package `ILAMM', PWLAD proposed by \citet{gao18},
and PIQ proposed by \citet{she22}.

To investigate the numerical performance of our method, the following measures are considered over 100 simulations:

1. M: the masking probability (fraction of undetected true outliers),

2. S: the swamping probability (fraction of good points labeled as outliers),

3. JD: the joint outlier detection rate (fraction of simulations with 0 masking),

4. FZR: the false zero rate (fraction of nonzero coefficients that are estimated as zero),
$$
\mathrm{FZR}(\widehat{\boldsymbol{\beta}})=\left|\left\{j \in\{1, \ldots, p_n\}: \widehat{\beta}_{j}=0 \wedge \beta_{j} \neq 0\right\}\right| /\left|\left\{j \in\{1, \ldots, p_n\}: \beta_{j} \neq 0\right\}\right|,
$$
where $|S|$ denotes the size of set $S$.

5. FPR: the false positive rate (fraction of zero coefficients that are estimated as nonzero),
$$
\operatorname{FPR}(\widehat{\boldsymbol{\beta}})=\left|\left\{j \in\{1, \ldots, p_n\}: \hat{\beta}_{j} \neq 0 \wedge \beta_{j}=0\right\}\right| /\left|\left\{j \in\{1, \ldots, p_n\}: \beta_{j}=0\right\}\right|.
$$

6. SR: the correct selection rate (fraction of identifying both nonzeros and zeros of $\boldsymbol{\beta}$),

7. CR: the correct coverage rate (fraction of identifying nonzeros of $\boldsymbol{\beta}$),


8. EE: the sum of squares of estimate error 
$$
\mathrm{EE}=\left(\widehat{\boldsymbol{\beta}}-\boldsymbol{\beta}_0\right)^{\mathrm{T}} \left(\widehat{\boldsymbol{\beta}}-\boldsymbol{\beta}_0\right).
$$

9. EE$_{non}$: the sum of squares of estimate error restrict to the important variables (nonzero coefficients)

For better performance in terms of outlier detection, it is preferable to maintain minimal values for both M and S, while simultaneously maximizing JD. In the context of sparse estimator $\boldsymbol{\beta}$, superior performance in variable selection is indicated by lower FZR and FPR values, as well as elevated higher SR and CR values. Concerning estimation accuracy, it is ideal for EE and EE$_{non}$ to be as minimal as possible.

\subsection{Simulation studies}
In this section, we consider both homogeneous model and heterogeneous model. We simulated data from the high dimensional model, accounting for the presence of outliers in response variables, covariates, and their combinations 
We set $n=100$ and $p_n=400$. The true coefficient is set as $\beta=(0.8,0.8,0.8,\ldots,0)^{\prime}$ with the number of nonzero components $d_n=3$, 
and the remaining ($p_n-d_n$) elements being zero. 

\Rmnum{1} (Homogeneous model). We generate covariate matrix $X=\left(x_{1}, \ldots, x_{n}\right)^{\prime}$ from a multivariate normal distribution $N(0,\boldsymbol{\Sigma})$, where
$\boldsymbol{\Sigma}=\left(\Sigma_{k m}\right)_{p_n \times p_n}$ with $\Sigma_{k m}=0.5^{|k-m|}$.
Continuous response variables are generated in accordance with model (\ref{linear}) for $i=1,\cdots,n$. 
The error term $\varepsilon_i$ is independent with $x_i$, and i.i.d according to the following distributions:

T1. the standard normal distribution $N(0,1)$.

T2. the Students' t-distribution with three degrees of freedom ($t_3$). 

Without loss of generality, we only contaminate the first $[cn]$ ($0<c<1$) observations by following scenarios, producing contaminated data set $(y^{*}_{i},x^{*}_{i})$ for $i=1,\cdots,n$.
We consider $c=0.1, 0.2, 0.3$ for all settings.

Case 1. The outlying perturbation was introduced on the response. The influential observations are generated according to
${y}^{*}_{i}=x^{*\top}_{i} \widetilde{\boldsymbol{\beta}}+\varepsilon_{i}$ for $i=1,\ldots,[cn]$, where 
$\widetilde{\boldsymbol{\beta}}=(0.8,0.8,0.8, \kappa, \kappa, \ldots, \kappa)^{\top}$
and $x^{*}_{i}=x_{i}$. This equivalent to add influential points $\kappa Z_{i}$ with $Z_{i}=x_{i}^{\top} \gamma$ and $\gamma=(0,0,0,1,1, \ldots, 1)^{\top}$ on the 
$i$-th observation. 

Case 2. We introduce the perturbation on covariates while keeping the response uncontaminated. 
We set $x^{*}_{i j}=x_{i j}+30 \kappa I_{\{1\leq j \leq 100\}}$ for the first $[cn]$ observations.

Case 3. The perturbation was introduced on both the response and covariates, that is $(y^{*}_{i},x^{*}_{ij})$ for
$i=1,\ldots,[cn], \quad j=1,\ldots,100$ under same contamination in Case 1 and Case 2.

The parameter $\kappa$ determines the impact of influential points, and here we set $\kappa=0.4$. 
The results are given in Table \ref{tab1}($c=0.1$) and Tables D1-D2($c=0.2,0.3$) in the Appendix D.

From Table \ref{tab1}, it is evident that PM exhibits a more pronounced masking effect in outlier detection compared to others when leverage points are present in the data, while PWLAD experiencing a substantial swamping effect.  PIQ performs well when there are high leverage points and combinations, however, its detect on strong masking when the responses
	of the influential observations are contaminated by a random perturbation.

Regarding robust parameter estimation and variable selection, PM displays notable limitations when leverage points are present and/or error terms adhere to a heavy-tailed distribution, 
characterized by lower SR, CR, and higher FZR. PWLAD's performance deteriorates when errors possess a heavy tail, leading to increased estimate error. 
Meanwhile, the AHuber estimators' efficiency in variable selection diminishes, evidenced by a reduced CR and elevated FZR when outliers are present in both response and covariates.
 For PIQ, it loses efficiency on variable selection, although it is conducted using the same threshold as PWHL.

As the contamination proportion progressively rises (see Table D1 and D2 in the Appendix D), PM's efficiency in both outlier detection and variable selection declines, despite its high breakdown point being approximately optimal to $\frac{1}{2}$. Conversely, the proposed estimator maintains robustness under conditions of high contamination and heavy-tailed noise scenarios.

\Rmnum{2} (Heterogeneous model). We generate continuous response variable according to the following model:
\begin{equation}
	y_{i}=x_{i}^{\top}\boldsymbol{\beta} +\eta_{i} \varepsilon_{i},
\end{equation}
where $\eta_i=\exp(x_{i20}+x_{i21}+x_{i22})$ generate schemes including heteroscedasticity. 
Other settings are the same as Simulation \Rmnum{1}. 

The results are presented in Table \ref{tab2}, illustrating that PM, as anticipated, exhibits inferior performance in outlier detection, characterized by a lower JD and a higher masking effect than PWLAD and PWHL. In comparison to PWLAD, where outlier indicator weights are updated based on absolute deviation, the residuals stemming from the adaptive Huber's loss function partially mitigate the swamping effect. Consequently, PWHL results display superiority, exhibiting a lower S than the other two methods.

In terms of estimation accuracy, the PWLAD estimator yields a higher EE than PM and PWHL, as LAD loses estimation efficiency under asymmetric noise distribution. Employing appropriate data-driven tuning parameters, the proposed estimator not only achieves the lowest estimate error (EE and EE$_{non}$) but also the highest correct variable selection rate (SR and CR) among the five robust methods.

\begin{footnotesize}
	\begin{longtable}{lcccccccccc}
		\caption{Comparison of PM\citep{kong18}, PWLAD, AHuber, PIQ and the proposed method PWHL
			under three cases in homogeneous model for $c=0.1$.}\label{tab1}\\
		\hline
		(T, Case)  & Method & M & S & JD &  EE & EE$_{non}$  & FZR & FPR & SR & CR \\
		\hline
		\endfirsthead	
		\hline
		(T, Case)  & Method & M & S & JD &  EE & EE$_{non}$  & FZR & FPR & SR & CR \\
		\hline
		\endhead
		\hline
		\endfoot
		(T1, Case 1) & PM & 0.2220 & \textbf{0.0033} & 0.11 & 0.4922 & 12.0265 & 0.0233 & 0.0008 & 0.67 & 0.94 \\ 
		& AHuber & - & - & -& 0.6276 & 11.5431 & 0.0500 & 0.0177 & 0.17 & 0.89 \\ 
		& PWLAD & \textbf{0.0200} & 0.4667 & \textbf{0.80} & 0.5731 & 13.8329 & 0.0000 & 0.0635 & 0.00 & 1.00 \\ 
		& PWHL & {0.0900} & 0.1789 & 0.40 & \textbf{0.4730} & \textbf{0.4730} & 0.0000 & 0.0000 & \textbf{1.00} & \textbf{1.00} \\
		&PIQ & 0.1900 & 0.0767 & 0.00 & 0.7123 & 0.6163 & 0.2000 & 0.0040 & 0.00 & 0.40 \\
		\hline
		(T1, Case 2) & PM & 0.9850 & 0.0001 & 0.00 & \textbf{0.2829} & 14.3199 & 0.0000 & 0.0100 & 0.03 & 1.00 \\ 
		& AHuber & - & - & - & 1.2816 & 2.0190 & 0.0000 & 0.2983 & 0.00 & 1.00 \\ 
		& PWLAD & 0.5980 & 0.4071 & 0.14 & 0.5436 & 13.8971 & 0.0000 & 0.0762 & 0.00 & 1.00 \\  
		& PWHL & \textbf{0.0000} & \textbf{0.0004} & \textbf{1.00} & 0.5204 & 0.5672 & {0.0000} & \textbf{0.0002} & \textbf{0.94} & 1.00 \\
		& PIQ & {0.0000} & 0.0556 & {1.00} & 0.4803 & \textbf{0.3985} & 0.0667 & 0.0030 & 0.00 & 0.80 \\
		\hline
		(T1, Case 3) & PM & 0.0000 & 0.0000 & 1.00 & 0.3284 & 13.4389 & 0.0033 & 0.0004 & 0.84 & 0.99 \\
		& AHuber & - & - & - & 1.3856 & 1.3856 & 1.0000 & 0.0000 & 0.00 & 0.00 \\ 
		& PWLAD & 0.0000 & 0.2544 & 1.00 & 0.5834 & 13.7694 & 0.0000 & 0.0776 & 0.00 & 1.00 \\
		& PWHL & 0.0000 & \textbf{0.0022} & 1.00 & \textbf{0.3879} & \textbf{0.3879} & 0.0000 & 0.0000 & \textbf{1.00} & 1.00 \\
		& PIQ & 0.0000 & 0.0556 & 1.00 & 7.1602 & 1.4331 & 0.9833 & 0.0099 & 0.00 & 0.00 \\
		\hline
		(T2, Case 1) & PM & 0.2890 & \textbf{0.0230} & 0.12 & \textbf{0.6663} & 10.6402 & 0.0933 & 0.0006 & 0.60 & 0.80 \\
		& AHuber & - & - & - & 0.8096 & 10.2474 & 0.1367 & 0.0152 & 0.14 & 0.74 \\ 
		& PWLAD & \textbf{0.0900} & 0.3367 & 0.30 & 0.8269 & 12.3654 & \textbf{0.0000} & 0.0763 & 0.00 & \textbf{1.00} \\
		& PWHL & 0.1100 & 0.1956 & 0.30 & 0.8147 & 0.8139 & 0.0667 & \textbf{0.0005} & \textbf{0.60} & 0.80 \\
		&PIQ & 0.1680 & 0.0742 & 0.00 & 0.9579 & \textbf{0.7976} & 0.3467 & 0.0051 & 0.00 & 0.31 \\
		\hline
		(T2, Case 2) & PM & 0.9570 & 0.0098 & 0.00 & \textbf{0.4532} & 13.6776 & 0.0033 & 0.0098 & 0.02 & 0.99 \\
		& AHuber & - & - & - & 1.2084 & 2.9423 & 0.0000 & 0.3005 & 0.00 & 1.00 \\ 
		& PWLAD & 0.4800 & 0.5017 & 0.23 & 0.6779 & 13.7085 & 0.0000 & 0.0759 & 0.00 & 1.00 \\
		& PWHL & \textbf{0.0000} & \textbf{0.0067} & \textbf{1.00} & 0.7281 & 1.0347  & \textbf{0.0000} & \textbf{0.0013} & \textbf{0.60} & 1.00 \\ 
		& PIQ & 0.0000 & 0.0556 & 1.00 & 0.8276 & \textbf{0.6899} & 0.2600 & 0.0045 & 0.00 & 0.39 \\
		\hline
		(T2, Case 3) & PM & 0.0000 & 0.0150 & 1.00 & \textbf{0.5153} & 12.1568 & 0.0233 & 0.0005 & 0.79 & 0.95 \\ 
		& AHuber & - & - & -& 1.3856 & 1.3856 & 1.0000 & \textbf{0.0000} & 0.00 & 0.00 \\ 
		& PWLAD & 0.0000 & 0.3511 & 1.00 & 0.6902 & 13.6898 & 0.0000 & 0.0695 & 0.00 & 1.00 \\   
		& PWHL & 0.0000 & \textbf{0.0044} & 1.00 & 0.5906 & \textbf{0.5905} & \textbf{0.0000} & 0.0005 & \textbf{0.80} & 1.00 \\
		& PIQ & 0.0000 & 0.0556 & 1.00 & 7.2011 & 1.4910 & 0.9867 & 0.0100 & 0.00 & 0.00 \\
		\hline
	\end{longtable}	
\end{footnotesize}

\begin{footnotesize}
	\begin{longtable}{lcccccccccc}
		\caption{Comparison of PM\citep{kong18}, PWLAD, AHuber, PIQ and the proposed method PWHL
			under three cases in heterogeneous model for $c=0.1$.}\label{tab2}\\
		\hline
		(T, Case)  & Method & M & S & JD &  EE & EE$_{non}$  & FZR & FPR & SR & CR \\
		\hline
		\endfirsthead	
		\hline
		(T, Case)  & Method & M & S & JD &  EE & EE$_{non}$  & FZR & FPR & SR & CR \\
		\hline
		\endhead
		\hline
		\endfoot
		(T1, Case 1) & PM & 0.2310 & 0.2210 & 0.15 & 0.7455 & 9.9876 & 0.2333 & 0.0013 & 0.32 & 0.59 \\ 
		& AHuber & - & - &- & 1.8771 & 5.5841 & 0.6967 & 0.0344 & 0.01 & 0.10 \\  
		& PWLAD & \textbf{0.0900} & 0.4944 & \textbf{0.40} & 0.8733 & 12.5224 & \textbf{0.0667} & \textbf{0.0678} & 0.00 & \textbf{0.80} \\ 
		& PWHL & 0.1000 & 0.3911 & 0.20 & \textbf{0.7196} & \textbf{0.7196} & 0.1667 & 0.0000 & \textbf{0.70} & 0.70 \\
		& PIQ & 0.5300& \textbf{0.1144} &0.00 & 1.7678 & 1.0581 &0.4000&0.0055&0.00&0.20\\
		\hline
		(T1, Case 2) & PM & 0.4570 & 0.1906 & 0.09 & 0.6433 & 11.9388 & 0.0800 & 0.0070 & 0.10 & 0.72 \\ 
		& AHuber & - & - & - & 1.1282 & 7.0040 & 0.1000 & 0.1824 & 0.00 & 0.70 \\ 
		& PWLAD & 0.0800 & 0.8954 & 0.64 & 2.2912 & 6.5301 & 0.5667 & 0.0400 & 0.00 & 0.21 \\  
		& PWHL & \textbf{0.0000} & 0.1300 & \textbf{1.00} & \textbf{0.5038} & \textbf{0.4973} & \textbf{0.0000} & \textbf{0.0028} & \textbf{0.50} & \textbf{1.00} \\
		& PIQ & 0.0400& \textbf{0.0600} &0.00 & 3.0349 & 1.1999 &0.7000&0.0078&0.00&0.10\\
		\hline
		(T1, Case 3) & PM & 0.0000 & 0.2276 & 1.00 & 0.7139 & 10.1015 & 0.2233 & 0.0010 & 0.40 & 0.63 \\
		& AHuber & - & - & - & 1.8441 & 3.1425 & 0.9033 & 0.0157 & 0.00 & 0.09 \\
		& PWLAD & 0.0000 & 0.4356 & 1.00 & 0.7626 & 13.4110 & 0.0667 & 0.0673 & 0.00 & 0.80 \\ 
		& PWHL & 0.0000 & 0.1344 & 1.00 & \textbf{0.4732} & \textbf{0.4732} & \textbf{0.0333} & \textbf{0.0000} & \textbf{0.90} & \textbf{0.90} \\
		& PIQ & 0.0000& \textbf{0.0056} &1.00 & 7.6911 & 1.3856 &1.0000&0.0101&0.00&0.00\\
		\hline
		(T2, Case 1) & PM & 0.3280 & 0.2177 & 0.09 & 0.9026 & 8.5582 & 0.3667 & 0.0008 & 0.26 & 0.43 \\ 
		& AHuber & - & - & - & 2.3409 & 5.7085 & 0.7533 & 0.0414 & 0.00 & 0.05 \\
		& PWLAD & 0.0900 & 0.5556 & \textbf{0.50} & 0.8911 & 12.6949 & \textbf{0.0333} & 0.0698 & 0.00 & \textbf{0.90} \\
		& PWHL & \textbf{0.0700} & 0.3900 & 0.40 & \textbf{0.7542} & \textbf{0.7524} & 0.1000 & \textbf{0.0008} & \textbf{0.70} & 0.80 \\ 
		& PIQ & 0.5850& \textbf{0.1206} &0.00 & 2.0357 & 1.1851 &0.5433&0.0066&0.00&0.09\\
		\hline
		(T2, Case 2) & PM & 0.5390 & 0.2321 & 0.14 & 0.6508 & 12.4313 & 0.0933 & 0.0073 & 0.12 & 0.78 \\ 
		& AHuber & - & - & - & 1.1241 & 5.7844 & 0.2033 & 0.1393 & 0.01 & 0.75 \\ 
		& PWLAD & 0.0520 & 0.9006 & 0.81 & 2.9619 & 5.2771 & 0.6700 & 0.0373 & 0.00 & 0.17 \\ 
		& PWHL & \textbf{0.0000} & 0.1478 & \textbf{1.00} & \textbf{0.5625} & \textbf{0.5605} & \textbf{0.0000} & \textbf{0.0013} & \textbf{0.50} & \textbf{1.00} \\
		& PIQ & 0.0660& \textbf{0.0629} &0.00 & 3.8435 & 1.4565 &0.7300&0.0080&0.00&0.01\\
		\hline
		(T2, Case 3) & PM & 0.0000 & 0.2280 & 1.00 & 0.8732 & 8.9798 & 0.3600 & \textbf{0.0008} & 0.27 & 0.45 \\ 
		& AHuber & - & - & -& 2.0017 & 4.2064 & 0.8967 & 0.0163 & 0.00 & 0.09\\
		& PWLAD & 0.0000 & 0.3489 & 1.00 & 0.9163 & 12.5154 & \textbf{0.0667} & 0.0632 & 0.00 & \textbf{0.80} \\
		& PWHL & 0.0000 & \textbf{0.1467} & 1.00 & \textbf{0.8394} & \textbf{0.8394} & 0.1667 & 0.0010 & \textbf{0.50} & 0.60 \\
		& PIQ & 0.0000 & 0.2000 & 1.00& 7.2773 & 1.3854 & 0.9967 & 0.0101 & 0.00 & 0.00 \\
		\hline
	\end{longtable}	
\end{footnotesize}

\subsection{Real data analysis}
We employ our proposed method to examine breast cancer survival data derived from the breast cancer study conducted by\citep{van06}, 
and illustrate the performances of PM, AHuber, PWLAD, PIQ, and PWHL. The study encompasses 295 female patients diagnosed with primary invasive breast carcinoma, 
with the expression of 24,885 genes profiled on cDNA arrays from all tumors. \citet{veer02} selected a set of 4,919 candidate genes through an initial screening utilizing the Rosetta error model.

Prior research has been conducted to identify genes associated with the overall survival of breast cancer patients, as evidenced in studies by\citet{song13} and \citet{zhang21}, where 
various screening methods have been proposed for high-dimensional survival data.
In this paper, we focus on selecting important genes related to the tumor diameter of breast cancer patients.

\begin{figure}
	\begin{center}
		\includegraphics[width=.7\textwidth]{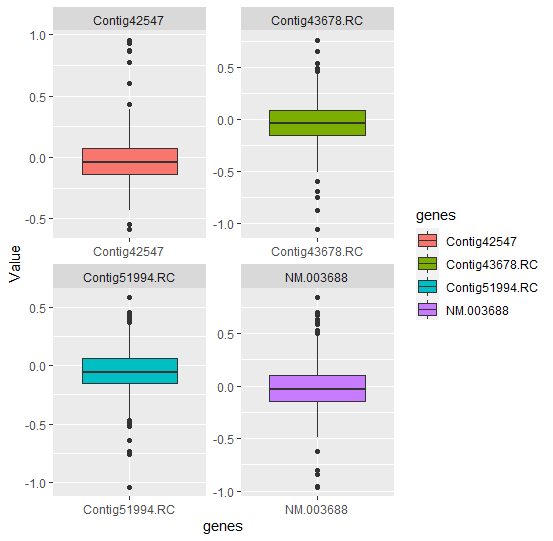}
	\end{center}
	\caption{The boxplots of four genes: Contig42547, Contig43678.RC, Contig51994.RC, and NM.003688.\label{fig1}}
\end{figure}

High-dimensional gene expression profiles frequently exhibit outliers, as depicted in Figure \ref{fig1}. \citet{song13} also identified potential outliers in gene expressions, which can result in biased estimation and prediction. To improve performance, \citet{song13} analyzed the data with outliers removed, an approach that may reduce the efficiency of estimators when a significant amount of data is discarded. Consequently, we reevaluate this dataset employing the previously mentioned robust methods to simultaneously select critical genes and detect outliers.

Initially, we screen the genetic expression dataset based on marginal association with tumor diameter for each gene.

Utilizing sure independence screening (SIS), a conservative selection of 590 variables is made, which can be executed using the "SIS" package in R to guarantee that all pertinent covariates pass this stage of the procedure. Subsequently, the remaining 590 covariates undergo further examination employing the aforementioned methods, yielding different gene reference sequences. 
Ultimately, we model tumor diameter in accordance with the linear transformation model:
\begin{equation}\label{realmodel}
	\log D_{i}=x_{i}^{\top} \boldsymbol{\beta}, \text { for } i=1, \ldots, n,
\end{equation}
where $D_i$ is the tumor diameter of the $i$-th subject, $x_i=(x_{i1},x_{i2},\ldots,x_{i590})$ is a 
590-dimensional genes vector with coefficient vector $\boldsymbol{\beta}=(\beta_1,\ldots,\beta_{590})^{\rm T}$.

We implement the five aforementioned methods to execute the regression. Residual inspections from PM and PWHL are displayed in Figures D1--D2 in Appendix D.2. The selected genes and their corresponding coefficients by the five methods are reported in Table \ref{tab3}(due to the limited space, we omit genes selected by PIQ, and more details are given in the Appendix D.2.). From the provided figures, the Residuals-Fitted plots indicate that the linear regression assumption for (\ref{realmodel}) is reasonable. The QQ-plots in both figures identify residuals from the 35th and 92nd observations of the fitted regressions, which may exhibit skewed distributions compared to the majority.

 For a further inspection of endowed heterogeneity, we plot Scale-Location graphs in Figures D1--D2 presented in Appendix D.2, and skewed line in Figure D1 suggests potential heteroscedasticity.   
Drawing from the simulation studies in Section 4, our proposed estimator tends to outperform others in addressing both heteroscedasticity and outliers.

As observed in Table \ref{tab3}, most genes selected by PM are also selected by the AHuber, PWHL, and PIQ methods, suggesting that these genes may have a significant impact on patients' tumor diameter.  
There is few overlap in the genes selected by PWLAD and the other four methods, except two genes(namely Contig50981.RC and Contig19284.RC) being selected by PIQ, indicating that the genes identified by PWLAD warrant further verification.

Table \ref{tab3} shows that PM selects 12 genes, PWLAD selects 14 genes, AHuber selects 42 genes, PIQ selects 171 genes, while our proposed PWHL selects 27 genes.
In terms of outlier detection, 
PM fails to detect any outliers, whereas PIQ considers all data as contamination, which contradicts the facts and demonstrates that PM and PIQ loses efficiency in outlier detection.
In comparison, PWLAD and our proposed procedure identify 202 and 181 potential outliers, respectively, with the 113th, 178th, 257th, and 235th observations detected as outliers.
This is consistent with our intuitive conclusions drawn from Figures D1--D2.

Simulation results in Section 4 indicate that when there are potential leverage points and/or residuals from the fitted regressions 
exhibit a heavy tail and skewed distribution, PWLAD experiences a more severe swamping effect than PWHL, which implies that outlier detection results using 
our proposed procedure can be more reliable. 

\begin{figure}[h]
	\begin{center}
		\includegraphics[width=.7\textwidth]{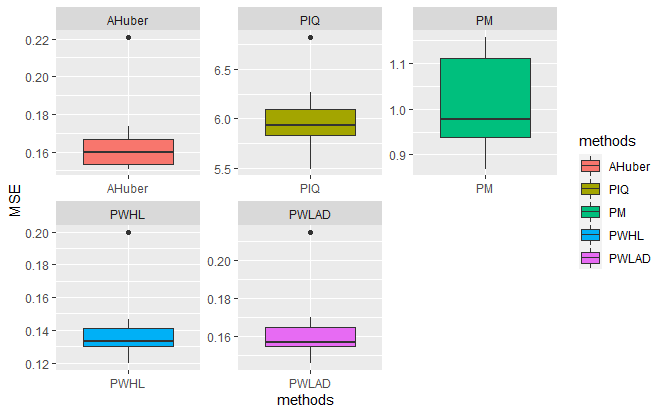}
	\end{center}
	\caption{Boxplots of the Mean Squared Error of predictions.\label{fig4}}
\end{figure}

To evaluate the predictive accuracy of PWHL, we employ the training set to construct the predictive model and use the testing data to evaluate the model.
We repeated the random splitting process 100 times, with 75\% of observations randomly selected as training data and the remainder allocated for testing.
Four methods mentioned earlier are applied to the training dataset to obtain the corresponding estimated coefficients, which are then used to predict the responses of the testing data. 
The boxplots of the Mean Squared Error (MSE) of predictions are illustrated in Figure \ref{fig4}, demonstrating that the proposed method exhibits the lowest MSE among all competitors. All these analysis results indicate that our proposed PWHL method performs commendably in robust ultrahigh-dimensional variable selection and offers favorable prediction performance.

\begin{footnotesize}
	\setlength{\tabcolsep}{1pt}
	\begin{longtable}{lcccccc}
		\caption{Selected genes and corresponding coefficients by different methods}\label{tab3}\\
		\hline
		\multicolumn{1}{c}{ Method}&
		\multicolumn{6}{c}{ Selected genes}\\
		\hline
		\endfirsthead	
		\hline
		\multicolumn{1}{c}{ Method}&
		\multicolumn{6}{c}{ Selected genes}\\
		\hline
		\endhead
		\hline
		\endfoot
		PM & NM.006115 & NM.000668 & NM.002989 & Contig14284.RC & NM.017680 & X99142 \\
		& 0.1478& -0.0355 & -0.2286 & -0.0180 & -0.0629 & 0.1153 \\
		& NM.000111 & Contig24609.RC & NM.001267 & NM.001395 & Contig58301.RC & NM.005101 \\
		& -0.0654 & -0.0321 & -0.0497 & -0.0399 & 0.0195 & 0.0305 \\ 
		\hline
		PWLAD & NM.000291 & Contig50981.RC & NM.007274 & AL137718 & NM.017450 & D15050 \\
		& 0.0187 & 0.1627 & 0.2040 & -0.0063 & 0.3596 & 0.2831 \\ 
		& AL080079 & NM.003688 & NM.007268 & NM.006931 & NM.013994 & NM.014736 \\
		& 0.1346 & 0.3888 & 0.2378 & 0.2216 & 0.3479 & 0.1361 \\ 
		& Contig19284.RC & NM.002509 & & & & \\
		& 0.5886 & 0.0492 &     &    & & \\
		\hline
		AHuber & NM.006115 & NM.003258 & NM.004867 & NM.003022 & NM.001168 & NM.002989 \\
		& 0.1116 & 0.4546 & -0.0522 & -0.1153  & -0.2076  & -0.1792 \\ 
		& NM.014321 & NM.006829 & AL133619 & NM.001809 & NM.018407 & Contig14284.RC \\
		& 0.0554 & -0.1643 & -0.0705 & -0.1835 & -0.0527 & 0.0544 \\ 
		& NM.000237 & NM.006096 & NM.017680 & Contig45537.RC & NM.005410 & AL157502 \\
		& -0.1709 & 0.0599 & -0.0122 & -0.0549 & 0.0453  & 0.1389 \\ 
		& Contig44191.RC & Contig53307.RC & X99142 & Contig57138.RC & NM.000111 & NM.004780 \\
		& 0.1420 & 0.0470 & 0.0711 & -0.0134 & -0.0517 & 0.1600 \\ 
		& Contig7558.RC & NM.006623 & NM.000662 & AF079529 & AI147042.RC & Contig46597.RC \\
		& -0.1336 & -0.0190 & 0.0208 & 0.0003 & 0.1014 & 0.0503 \\ 
		& Contig24609.RC & K02403 & Contig20629.RC & NM.001394 & Contig41850.RC & NM.007088 \\
		& -0.0209 & -0.0113 & -0.1001 & -0.0114 & -0.0384 & 0.0702 \\ 
		& NM.001267 & X03084 & NM.001395 & Contig58301.RC & NM.005101 & NM.014373 \\
		& -0.0293 & 0.1274 & 0.0052 & 0.0290 & 0.0205  & -0.0868 \\  
		\hline
		PWHL & NM.006115 & NM.003022 & NM.001168 & NM.000668 & NM.002989 & NM.006829 \\
		& 0.1070 & -0.1134 & 0.0274 & -0.0522 & -0.2247 & -0.0894 \\ 
		& AL133619 & Contig14284.RC & NM.006096 & NM.017680 & NM.005410 & AL157502 \\
		& -0.0847 & 0.0524 & 0.0518 & -0.0113 & 0.0362 & 0.0748 \\ 
		& Contig44191.RC & X99142 & Contig57138.RC & NM.000111  & AF079529 & AI147042.RC \\
		& 0.1114 & 0.0796 & 0.0000 & -0.0690 & -0.0372 & 0.0978 \\ 
		& Contig24609.RC & NM.001394 & Contig41850.RC & NM.007088 & NM.001267 & NM.001395\\
		& -0.0312 & -0.0033 & -0.1157 & 0.0126 & -0.0296 & -0.0540 \\ 
		& Contig58301.RC & NM.005101 & NM.014373 &     &    &  \\
		& 0.0433 & 0.0299 & -0.0664 &     &    &  \\ 	
		\hline
	\end{longtable}	
\end{footnotesize}

\section{Conclusions}
This article develops a penalized weighted Huber-LASSO procedure for robust variable selection and outliers detection.
Instead of estimating the conditional quantile regression as in PWLAD \citep{gao18}, we construct adaptive Huber's estimators based on the dependent mean regression function,
accommodating a heteroscedastic linear model with asymmetric noise distribution. As a hybrid of squared loss(LS) for relatively small errors and absolute loss (LAD) for rather large errors, 
the adaptive Huber's procedure governs the degree of hybridization using a single data-driven tuning parameter, facilitating ultrahigh-dimensional mean regression even in the presence of heavy-tailed errors.

Moreover, outliers in high-dimensional datasets may result in severe estimation bias, which is a notable limitation of least squares estimators. We compute robust estimators resistant to both outliers present in response variables and leverage points, offering superiority over some procedures based on the mean shift model (see \citet{kong18}). Regarding outlier detection, the weight vectors serve as an instructor for latent outliers with values less than 1. Additionally, the weights are updated using adaptively truncated residuals encapsulated within Huber's loss function, mitigating the swamping effect compared to using LS or LAD loss, particularly in heterogeneous data.

Our numerical studies demonstrate that the proposed procedure achieves high effectiveness in robust variable selection and is robust for both response variables and covariates in high-dimensional linear models. It performs exceptionally competitive under heavy-tailed distributions and/or heteroscedastic models, with promising prospects when the dimension of covariates is exponentially larger than the sample size.

Our proposed procedure can be extended to generalized linear models for further study, while extensive simulations would be necessary to assess performance.
The effects of initial values cannot be ignored, as they significantly impact the consistency of outlier detection. Nevertheless, residual-based outlier detection methods 
inherently exhibit masking and/or swamping effects when detecting multiple outliers in high dimensions. Accurately detecting multiple moderate outliers in high-dimensional settings remains a challenging problem.

\section*{Acknowledgments}
The authors thank the Associate Editor and three referees for their constructive comments.

\par


\bibhang=1.7pc
\bibsep=2pt
\fontsize{9}{14pt plus.8pt minus .6pt}\selectfont
\renewcommand\bibname{\large \bf References}
\expandafter\ifx\csname
natexlab\endcsname\relax\def\natexlab#1{#1}\fi
\expandafter\ifx\csname url\endcsname\relax
\def\url#1{\texttt{#1}}\fi
\expandafter\ifx\csname urlprefix\endcsname\relax\def\urlprefix{URL}\fi



\newpage

\appendix
\setcounter{table}{0}   
\setcounter{figure}{0}
\renewcommand{\thetable}{D\arabic{table}}
\renewcommand{\thefigure}{D\arabic{figure}}

\section{Proof of Theorems}

{\bf Proof of Theorem 1.} To simplify our proof, here we only consider homogeneous model, and the following techniques can be also applied to heterogeneous model by applying Law of total expectation.
Under an initial estimate $\widehat{\boldsymbol{\beta}}^{(0)}$, the weight $\boldsymbol{w}$ can be estimated by
\begin{equation}
	\widehat{\boldsymbol{w}}(\mu ; \alpha)
	= \underset{\boldsymbol{w}}{\arg \min }\left\{\frac{1}{2} \sum_{i=1}^{n} \ell_{\alpha}\left({w_{i}}\left(y_{i}-x_{i}^{\top} \widehat{\boldsymbol{\beta}}^{(0)}\right)\right)+\mu \sum_{i=1}^{n} {\varpi}_{i}\left|1-w_{i}\right|\right\},
\end{equation}
where the estimate $\widehat{\boldsymbol{w}}$ satisfies
\begin{equation}
	\left\{\begin{array}{ll}
		\widehat{w}_{i}\ell_{\alpha}\left(y_{i}-x_{i}^{\top} \widehat{\boldsymbol{\beta}}^{(0)}\right)=\mu\varpi_{i}  & \text { if } 0<\widehat{w}_{i}<1 \\
		
		\widehat{w}_{i}\ell_{\alpha}\left(y_{i}-x_{i}^{\top} \widehat{\boldsymbol{\beta}}^{(0)}\right)<\mu \varpi_{i} & \text { if } \widehat{w}_{i}=1
	\end{array}\right.,
\end{equation}
according to our algorithm in Section 2.1.
To prove the outlier detection consistency, first we notice that
$$
\begin{aligned}
	P(\widehat{\mathcal{O}} \neq \mathcal{O}) & \leq P\left(\bigcup_{i \in \mathcal{O}}\left\{\ell_{\alpha}\left(y_{i}-x_{i}^{\top} \widehat{\boldsymbol{\beta}}^{(0)}\right)<\mu \varpi_{i}\right\}\right) \\
	&+P\left(\bigcup_{i \in \mathcal{O}^{c}}\left\{\ell_{\alpha}\left(y_{i}-x_{i}^{\top} \widehat{\boldsymbol{\beta}}^{(0)}\right) \geq \mu \varpi_{i}\right\}\right)\\
	& \triangleq \Rmnum{1}+\Rmnum{2}
\end{aligned}.
$$ 
For $\Rmnum{1}$, we have
$$
\begin{aligned}
	\Rmnum{1}&=P\left(\left|\hat{r}_{i}\right|>\alpha^{-1}\right)\cdot P\left(\bigcup_{i \in \mathcal{O}}\left\{\ell_{\alpha}\left(y_{i}-x_{i}^{\top} \widehat{\boldsymbol{\beta}}^{(0)}\right)<\mu \varpi_{i}\right\}I_{\left\{\left|\hat{r}_{i}\right|>\alpha^{-1}\right\}}\right)\\
	&+P\left(\left|\hat{r}_{i}\right|<\alpha^{-1}\right)P\left(\bigcup_{i \in \mathcal{O}}\left\{\ell_{\alpha}\left(y_{i}-x_{i}^{\top} \widehat{\boldsymbol{\beta}}^{(0)}\right)<\mu \varpi_{i}\right\}I_{\left\{\left|\hat{r}_{i}\right|<\alpha^{-1}\right\}}\right)\\
	& \triangleq \Rmnum{1}_1+\Rmnum{1}_2
\end{aligned},
$$
where $\hat{r}_{i}=y_i-x_i^{\top}\widehat{\boldsymbol{\beta}}^{(0)}$.
Under Conditions C2 and C3(i), we have,
\begin{equation}\label{11}
	\begin{aligned}
		\Rmnum{1}_1&\leq P\left(\bigcup_{i \in \mathcal{O}}\left\{2 \alpha^{-1}\left|\hat{r}_{i}\right| -\alpha^{-2}<\mu \cdot \varpi_{i}\right\}\right)\\
		&\leq P\left(\bigcup_{i \in \mathcal{O}}\left\{\left|\hat{r}_{i}\right|<\frac{\alpha}{2}\mu \cdot \varpi_{i}\right\}\right)\\
		&\leq P\left(\bigcup_{i \in \mathcal{O}}\left\{\left|\varepsilon_i\right|<\frac{\alpha}{2}\mu \cdot \underline{\varpi}_{n}+|x^{\top}_i(\widehat{\boldsymbol{\beta}}^{(0)}-\boldsymbol{\beta}^{*})|\right\}\right)+P\left(\max _{i \in \mathcal{O}} \varpi_{i}>\underline{\varpi}_{n}\right)\\
		& \leq q_{n} P\left(\left|\varepsilon_{i}\right|< \alpha \mu \underline{\varpi}_{n} \right)+P\left(\left\|\widehat{\boldsymbol{\beta}}^{(0)}-\boldsymbol{\beta}^{*}\right\|_{2}>\frac{\alpha}{2} \mu\underline{\varpi}_{n} /\left(l \sqrt{p_n}\right)\right)+o(1)\\
		& \leq q_n \alpha\cdot\mu\cdot\underline{\varpi}_{n}+P\left(\left\|\widehat{\boldsymbol{\beta}}^{(0)}-\boldsymbol{\beta}^{*}\right\|_{2}>\frac{\alpha}{2} \mu\underline{\varpi}_{n} /\left(l \sqrt{p_n}\right)\right)\\
		& \leq P\left(\left\|\widehat{\boldsymbol{\beta}}^{(0)}-\boldsymbol{\beta}^{*}\right\|_{2}>\frac{\alpha}{2} \mu\underline{\varpi}_{n} /\left(l \sqrt{p_n}\right)\right)+o(1)
	\end{aligned}.
\end{equation} 
Similarly, for $\Rmnum{1}_2$, we have
$$
\begin{aligned}
	\Rmnum{1}_2&\leq q_n P\left(\left|\varepsilon_i\right|<\alpha^{-1}+|x^{\top}_i(\widehat{\boldsymbol{\beta}}^{(0)}-\boldsymbol{\beta}^{*})|\right)\\
	&\leq q_n P\left(\left|\varepsilon_i\right|<\alpha^{-1}+\frac{\alpha}{2}\mu\underline{\varpi}_{n}\right)+P\left(\left\|\widehat{\boldsymbol{\beta}}^{(0)}-\boldsymbol{\beta}^{*}\right\|_{2}>\frac{\alpha}{2} \mu\underline{\varpi}_{n} /\left(l \sqrt{p_n}\right)\right)\\
	&\leq q_n\left(\alpha^{-1}+\frac{\alpha}{2}\mu\underline{\varpi}_{n}\right)+P\left(\left\|\widehat{\boldsymbol{\beta}}^{(0)}-\boldsymbol{\beta}^{*}\right\|_{2}>\frac{\alpha}{2} \mu\underline{\varpi}_{n} /\left(l \sqrt{p_n}\right)\right).
\end{aligned}.
$$
Under Condition C3(i), we have 
\begin{equation}\label{12}
	\Rmnum{1}_2\leq P\left(\left\|\widehat{\boldsymbol{\beta}}^{(0)}-\boldsymbol{\beta}^{*}\right\|_{2}>\frac{\alpha}{2} \mu\underline{\varpi}_{n} /\left(l \sqrt{p_n}\right)\right)+o(1).
\end{equation}

Next, we analyze $\Rmnum{2}$.
Following foregoing decomposition method, we have 
$$
\begin{aligned}
	\Rmnum{2}&=P\left(\left|\hat{r}_{i}\right|>\alpha^{-1}\right)\cdot P\left(\bigcup_{i \in \mathcal{O}^c}\left\{\ell_{\alpha}\left(y_{i}-x_{i}^{\top} \widehat{\boldsymbol{\beta}}^{(0)}\right)\geq\mu \varpi_{i}\right\}I_{\left\{\left|\hat{r}_{i}\right|>\alpha^{-1}\right\}}\right)\\
	&+P\left(\left|\hat{r}_{i}\right|<\alpha^{-1}\right)P\left(\bigcup_{i \in \mathcal{O}^c}\left\{\ell_{\alpha}\left(y_{i}-x_{i}^{\top} \widehat{\boldsymbol{\beta}}^{(0)}\right)\geq\mu \varpi_{i}\right\}I_{\left\{\left|\hat{r}_{i}\right|<\alpha^{-1}\right\}}\right)\\
	& \triangleq \Rmnum{2}_1+\Rmnum{2}_2
\end{aligned}.
$$
For $\Rmnum{2}_1$, we have
$$
\begin{aligned}
	\Rmnum{2}_1&\leq P\left(\bigcup_{i \in \mathcal{O}^c}\left\{2 \alpha^{-1}\left|\hat{r}_{i}\right| -\alpha^{-2}\geq \mu \cdot \varpi_{i}\right\}I_{\left\{\left|\hat{r}_{i}\right|>\alpha^{-1}\right\}}\right)\\
	& \leq  P\left(\bigcup_{i \in \mathcal{O}^c}\left\{2 \alpha^{-1}\left|\hat{r}_{i}\right|\geq \mu \cdot \varpi_{i}\right\}\right)\\
	& \leq P\left(\max _{i \in \mathcal{O}^{c}}\left|\varepsilon_{i}\right| \geq \frac{\alpha}{4}\mu \bar{\varpi}_{n}\right)+P\left(\left\|\widehat{\boldsymbol{\beta}}^{(0)}-\boldsymbol{\beta}^{*}\right\|_{2}>\frac{\alpha}{4} \mu \bar{\varpi}_{n}/\left(l \sqrt{p_n}\right)\right)\\
	& \leq \frac{\sqrt{2\log(2n)}\sigma}{\alpha/4\mu\bar{\varpi}_{n}}+P\left(\left\|\widehat{\boldsymbol{\beta}}^{(0)}-\boldsymbol{\beta}^{*}\right\|_{2}>\frac{\alpha}{4} \mu \bar{\varpi}_{n}/\left(l \sqrt{p_n}\right)\right)
\end{aligned}.
$$
Under Condition C3(ii), we have 
\begin{equation}\label{21}
	\Rmnum{2}_1 \leq P\left(\left\|\widehat{\boldsymbol{\beta}}^{(0)}-\boldsymbol{\beta}^{*}\right\|_{2}>\frac{\alpha}{4} \mu \bar{\varpi}_{n}/\left(l \sqrt{p_n}\right)\right).
\end{equation}
Finally, for $\Rmnum{2}_2$, there exists
\begin{equation}\label{22}
	\begin{aligned}
		\Rmnum{2}_2 &\leq P\left(\bigcup_{i \in \mathcal{O}^c}\left\{\alpha^{-1}|\hat{r}_i|\geq\mu \varpi_{i}\right\}\right)\\
		& \leq P\left(\max _{i \in \mathcal{O}^{c}}\left|\varepsilon_{i}\right| \geq \frac{\alpha}{2}\mu \bar{\varpi}_{n}\right)+P\left(\left\|\widehat{\boldsymbol{\beta}}^{(0)}-\boldsymbol{\beta}^{*}\right\|_{2}>\frac{\alpha}{4} \mu \bar{\varpi}_{n}/\left(l \sqrt{p_n}\right)\right)+o(1)\\
		& \leq \frac{\sqrt{2\log(2n)}\sigma}{\alpha/2\mu\bar{\varpi}_{n}}+P\left(\left\|\widehat{\boldsymbol{\beta}}^{(0)}-\boldsymbol{\beta}^{*}\right\|_{2}>\frac{\alpha}{4} \mu \bar{\varpi}_{n}/\left(l \sqrt{p_n}\right)\right)\\
		& \leq P\left(\left\|\widehat{\boldsymbol{\beta}}^{(0)}-\boldsymbol{\beta}^{*}\right\|_{2}>\frac{\alpha}{4} \mu \bar{\varpi}_{n}/\left(l \sqrt{p_n}\right)\right)+o(1)
	\end{aligned}.
\end{equation}
Combine initial estimate $P\left(\left\|\widehat{\boldsymbol{\beta}}^{(0)}-\boldsymbol{\beta}^{*}\right\|_{2}>\frac{\alpha}{2}\mu\max \left\{\bar{\varpi}_{n}/2, \underline{\varpi}_{n}\right\} /\left(\sqrt{p_n} l \right)\right)=o(1)$ with (\ref{11})--(\ref{22}), the theorem is proved.

{\bf Proof of Theorem 2.} Denote estimator $\widehat{\boldsymbol{\theta}}=\widehat{\boldsymbol{\theta}}_{{\mu},\lambda}$, and we introduce a new estimator defined as
\begin{equation}\label{1}
	\widehat{\boldsymbol{\theta}}_{\eta}=\boldsymbol{\theta}^{*}+\eta\left(\widehat{\boldsymbol{\theta}}-\boldsymbol{\theta}^{*}\right),
\end{equation}
with $\eta$ satisfying
$$
\eta=\left\{
\begin{array}{lll}
	1  & \text {if } & \left\|\Sigma^{1/2}\left(\widehat{\boldsymbol{\theta}}-\boldsymbol{\theta}^{*}\right)\right\|_{2} \leq r \\
	\in(0,1) & \text { otherwise } &\text { s.t. } \left\|\Sigma^{1/2}\left(\widehat{\boldsymbol{\theta}}_{\eta}-\boldsymbol{\theta}^{*}\right)\right\|_{2}=r
\end{array} \right.
$$
for some $r>0$.
We assert that when 
\begin{equation}\label{con}
	\left\|\nabla L_{\alpha}\left(\boldsymbol{\theta}^{*}\right)\right\|_{\infty}\leq \lambda/2,
\end{equation}
it holds that 
\begin{equation}\label{assert}
	\left\|\left(\widehat{\boldsymbol{\theta}}-\boldsymbol{\theta}^{*}\right)_{\mathcal{S}^{c}}\right\|_{1} \leq 3\left\|\left(\widehat{\boldsymbol{\theta}}-\boldsymbol{\theta}^{*}\right)_{\mathcal{S}}\right\|_{1}.
\end{equation}
For the optimal estimator $\widehat{\boldsymbol{\theta}}$, it holds that
\begin{equation}\label{1-1}
	L_{\alpha}(\widehat{\boldsymbol{\theta}})-L_{\alpha}\left(\boldsymbol{\theta}^{*}\right) \leq \lambda\left(\left\|\boldsymbol{\theta}^{*}\right\|_{1}-\|\widehat{\boldsymbol{\theta}}\|_{1}\right).
\end{equation}
In addition, it is noteworthy that
\begin{equation}\label{1-2}
	\begin{aligned}
		\|\widehat{\boldsymbol{\theta}}\|_{1}-\left\|\boldsymbol{\theta}^{*}\right\|_{1} 
		&=\left\|\widehat{\boldsymbol{\theta}}_{\mathcal{S}}+\widehat{\boldsymbol{\theta}}_{\mathcal{S}^{c}}\right\|_1-\left[\|\boldsymbol{\theta}_{S^{c}}^{*}\|_1+\|\boldsymbol{\theta}_{S}^{*}\|_1\right]\\ 
		&\geq\left\|\boldsymbol{\theta}_{S}^{*}\right\|_{1}+\left\|\left(\widehat{\boldsymbol{\theta}}-\boldsymbol{\theta}^{*}\right)_{\mathcal{S}^{c}}\right\|_1-\left\|\boldsymbol{\theta}_{S^{c}}^{*}\right\|_{1}-\left\|\left(\widehat{\boldsymbol{\theta}}-\boldsymbol{\theta}^{*}\right)_{\mathcal{S}}\right\|_{1}-\left(\left\|\boldsymbol{\theta}_{\mathcal{S}}^{*}\right\|_{1}+\left\|\boldsymbol{\theta}_{\mathcal{S}^{c}}^{*}\right\|_{1}\right) \\
		& \geq\left\|\left(\widehat{\boldsymbol{\theta}}-\boldsymbol{\theta}^{*}\right)_{\mathcal{S}^{c}}\right\|_{1}-\left\|\left(\widehat{\boldsymbol{\theta}}-\boldsymbol{\theta}^{*}\right)_{\mathcal{S}}\right\|_{1}
	\end{aligned}.
\end{equation}
Apply the first-order Taylor expansion to $L_{\alpha}(\widehat{\boldsymbol{\theta}})$ with regard to $\boldsymbol{\theta}^{*}$, combined with H$\ddot{o}$lder's Inequality, we have
\begin{equation}\label{1-3}
	\begin{aligned}
		L_{\alpha}(\widehat{\boldsymbol{\theta}})-L_{\alpha}\left(\boldsymbol{\theta}^{*}\right) & 
		\geq\nabla L_{\alpha}\left(\boldsymbol{\theta}^{*}\right)^{\top}\left(\widehat{\boldsymbol{\theta}}-\boldsymbol{\theta}^{*}\right)  \\
		&\geq-\left\|\nabla L_{\alpha}\left(\boldsymbol{\theta}^{*}\right)\right\|_{\infty}\|\widehat{\boldsymbol{\theta}}-\boldsymbol{\theta}^{*}\|_{1} \\
		& \geq-\frac{\lambda}{2}\left(\left\|\left(\widehat{\boldsymbol{\theta}}-\boldsymbol{\theta}^{*}\right)_{\mathcal{S}^{c}}\right\|_{1}+\left\|\left(\widehat{\boldsymbol{\theta}}-\boldsymbol{\theta}^{*}\right)_{\mathcal{S}}\right\|_{1}\right)
	\end{aligned}.
\end{equation}
Based on (\ref{1-1})--(\ref{1-3}), our assert (\ref{assert}) can be proved.
It is easy to prove that $\widehat{\boldsymbol{\theta}}_{\eta}$ also satisfy 
$$
\left\|\left(\widehat{\boldsymbol{\theta}}_{\eta}-\boldsymbol{\theta}^{*}\right)_{\mathcal{S}^{c}}\right\|_{1} \leq 3\left\|\left(\widehat{\boldsymbol{\theta}}_{\eta}-\boldsymbol{\theta}^{*}\right)_{\mathcal{S}}\right\|_{1},
$$ 
in accordance with convex optimization.
On the basis of Lemma C.4 in \citet{sun20}, under Conditions C6 and C7, it holds that
\begin{equation}\label{1-4}
	\begin{aligned}
		\left(\nabla L_{\alpha}\left(\widehat{\boldsymbol{\theta}}_{\eta}\right)-\nabla L_{\alpha}\left(\boldsymbol{\theta}^{*}\right)\right)^{\top}\left( \widehat{\boldsymbol{\theta}}_{\eta}-\boldsymbol{\theta}^{*}\right) 
		&\geq \frac{1}{4}\left\|\mathbf{\Sigma}^{1/2}\left(\widehat{\boldsymbol{\theta}}_{\eta}-\boldsymbol{\theta}^{*}\right)\right\|_{2}^{2} \\
		& \geq \frac{1}{4} \kappa_{l}^{1 / 2}\left\|\widehat{\boldsymbol{\theta}}_{\eta}-\boldsymbol{\theta}^{*}\right\|_{2}\left\|\mathbf{\Sigma}^{1/2}\left(\widehat{\boldsymbol{\theta}}_{\eta}-\boldsymbol{\theta}^{*}\right)\right\|_2
	\end{aligned}.
\end{equation}
with probability at least $1-(p_n+n)^{-1}$.
Also, it holds that
\begin{equation}\label{1-5}
	\left(\nabla \mathcal{L}_{\alpha}\left(\widehat{\boldsymbol{\theta}}_{\eta}\right)-\nabla \mathcal{L}_{\alpha}\left(\boldsymbol{\theta}^{*}\right)\right)^{\top} \left(\widehat{\boldsymbol{\theta}}_{\eta}-\boldsymbol{\theta}^{*}\right) 
	\leq \eta\left(\nabla \mathcal{L}_{\alpha}(\widehat{\boldsymbol{\theta}})-\nabla \mathcal{L}_{\alpha}\left(\boldsymbol{\theta}^{*}\right)\right)^{\top} \left(\widehat{\boldsymbol{\theta}}-\boldsymbol{\theta}^{*}\right),
\end{equation}
derived from Lemma C.1 in \citet{sun20}.
Moreover, from (\ref{1-1}) with (\ref{1-3}), it holds that
\begin{equation}
	\begin{aligned}
		&\left(\nabla L_{\alpha}(\widehat{\boldsymbol{\theta}})-\nabla L_{\alpha}\left(\boldsymbol{\theta}^{*}\right)\right)^{\top}\left(\widehat{\boldsymbol{\theta}}-\boldsymbol{\theta}^{*}\right)\\ 
		&\leq \lambda\left(\left\|\boldsymbol{\theta}^{*}\right\|_{1}-\|\widehat{\boldsymbol{\theta}}\|_{1}\right)+\frac{\lambda}{2}\left\|\widehat{\boldsymbol{\theta}}-\boldsymbol{\theta}^{*}\right\|_{1} \\
		&\leq\lambda\left(\left\|\left(\widehat{\boldsymbol{\theta}}-\boldsymbol{\theta}^{*}\right)_{\mathcal{S}}\right\|_{1}-\left\|\left(\widehat{\boldsymbol{\theta}}-\boldsymbol{\theta}^{*}\right)_{\mathcal{S}^{c}}\right\|_{1}\right)+\frac{\lambda}{2}\left\|\widehat{\boldsymbol{\theta}}-\boldsymbol{\theta}^{*}\right\|_{1}  \\
		&\leq\frac{\lambda}{2}\left(3\left\|\left(\widehat{\boldsymbol{\theta}}-\boldsymbol{\theta}^{*}\right)_{\mathcal{S}}\right\|_{1}-\left\|\left(\widehat{\boldsymbol{\theta}}-\boldsymbol{\theta}^{*}\right)_{\mathcal{S}^{c}}\right\|_{1}\right)
	\end{aligned}.
\end{equation}
Plug (\ref{1-4}) and (\ref{1-5}) into above result, yields
\begin{equation}
	\frac{1}{4} \kappa_{l}^{1 / 2}\left\|\widehat{\boldsymbol{\theta}}_{\eta}-\boldsymbol{\theta}^{*}\right\|_{2}\left\|\mathbf{\Sigma}^{1/2}\left(\widehat{\boldsymbol{\theta}}_{\eta}-\boldsymbol{\theta}^{*}\right)\right\|_{2} 
	\leq \frac{3}{2} \lambda s_n^{1 / 2}\left\|\eta \left(\widehat{\boldsymbol{\theta}}-\boldsymbol{\theta}^{*}\right)\right\|_{2}=\frac{3}{2} \lambda s_n^{1 / 2}\left\|\widehat{\boldsymbol{\theta}}_{\eta}-\boldsymbol{\theta}^{*}\right\|_{2},
\end{equation}
which implies that
\begin{equation}\label{l2}
	\left\|\mathbf{\Sigma}^{1/2}\left(\widehat{\boldsymbol{\theta}}_{\eta}-\boldsymbol{\theta}^{*}\right)\right\|_{2} \leq 6 \kappa_{l}^{-1 / 2} s_n^{1 / 2} \lambda
\end{equation}
for $l_2$-error upper bound, and
\begin{equation}\label{l1}
	\left\|\widehat{\boldsymbol{\theta}}_{\eta}-\boldsymbol{\theta}^{*}\right\|_{1} \leq 24 \kappa_{l}^{-1} s_n \lambda
\end{equation}
for $l_1$-norm oracle inequality.
According to Lemma C.6 in \citet{sun20}, when $\alpha^{-1}=\alpha^{-1}_{0}(n / \log(p_n+n))^{1 /(1+\delta)}$, it holds that
$$
P\left(\left\|\nabla L_{\alpha}\left(\boldsymbol{\theta}^{*}\right)\right\|_{\infty} \leq \left[(2\sqrt{2}+1)M+1\right] \sigma_{\max} \alpha^{-1}_{0}\left(\frac{\log(p_n+n)}{n}\right)^{\delta /(1+\delta)}\right)\leq 1-2(p_n+n)^{-1}.
$$
Therefore, when choosing $\lambda=2c \sigma_{\max} \alpha^{-1}_{0}\{(\log(p_n+n)) / n\}^{\delta /(1+\delta)}$ for some $c>\left[(2\sqrt{2}+1)M+1\right]$, and condition (\ref{con})
can be easily satisfied.
Furthermore, let $r=C(M,\kappa_{l}) \sigma_{\max} \alpha^{-1} \sqrt{(\log (p_n+n)) / n}$ with an appropriate constant $C(M,\kappa_{l})$ only depends on $M$ and $\kappa_{l}$, 
incorporate aforementioned conditions, it can be deduced that
$$
\left\|\mathbf{\Sigma}^{1/2}\left(\widehat{\boldsymbol{\theta}}_{\eta}-\boldsymbol{\theta}^{*}\right)\right\|_{ 2} \leq 12 c \kappa_{l}^{-1 / 2} \sigma_{\max} \alpha^{-1}_{0} s_n^{1 / 2}\left(\frac{\log (p_n+n)}{n}\right)^{\delta /(1+\delta)}<r
$$
with probability at least $1-3(p_n+n)^{-1}$. Let $\eta=1$ and we obtain the $l_2$-error upper bound for $\widehat{\boldsymbol{\theta}}$, and the $l_1$-norm oracle inequality (3.9) can be directly 
derived from (\ref{l1}). Applying Cauchy-Schwartz inequality, we obtain (3.10). Finally, we finish the proof.

{\bf Proof of Theorem 3.} We apply smoothed function technique proposed by \citet{heller07} to the estimating function , and the corresponding smoothed estimating function can be defined as:
\begin{equation}\label{sef}
	\begin{aligned}
		\widetilde{U}\left(Z;\boldsymbol{\theta}\right)&= \widetilde{\psi}_{\alpha}\left(Z, \boldsymbol{\theta}\right)+ p^{\prime}(\boldsymbol{\theta} ; \lambda)\\
		&=\left\{e\bar{\omega}(\alpha,h)+(1-\bar{\omega}(\alpha,h))\alpha^{-1}\left(2\Phi\left(\frac{e}{h}\right)-1\right)\right\}\Delta_{\theta}+p^{\prime}(\boldsymbol{\theta} ; \lambda)\\
		&\triangleq \widetilde{U}_F\left(Z;\boldsymbol{\theta}\right)+p^{\prime}(\boldsymbol{\theta} ; \lambda)
	\end{aligned}
\end{equation}
where $\bar{\omega}(\alpha, h)=\Phi\left(\frac{e+\alpha^{-1}}{h}\right)-\Phi\left(\frac{e-\alpha^{-1}}{h}\right)$. 

Since smoothing operation doesn't modify the convexity properties of Huber's loss function. According to Lemma 2 and Proposition 1 in \citet{marco17}, (1) and (2) in Theorem 3 can be proved
directly. To prove (3) and (4), we initially established a Lemma, which subsequentially allowed us to deduce the equivalence between 
the unsmoothed estimating function and the smoothed estimating equation (\ref{sef}).

{\bf Lemma. } Under Conditions C9--C11,
$$
\sup _{\theta}\left|U_F(Z;\boldsymbol{\theta})-\widetilde{U}_{F}(Z;\boldsymbol{\theta})\right| \stackrel{\text { a.s. }}{\rightarrow} 0.
$$

{\bf Proof of Lemma.} Let 
$$
\begin{aligned}
	g(e,Z)&=e\left(\bar{\omega}(\alpha, h)-I_{\{|e|\leq \alpha^{-1}\}}\right)+\alpha^{-1}\left(2\Phi(e/h)-1\right)\left(1-\bar{\omega}(\alpha, h)\right)\\
	&-\alpha^{-1}\sgn(e)I_{\{|e|> \alpha^{-1}\}}
\end{aligned}
$$
Since
\begin{equation}\label{inner}
	\begin{aligned}
		U_F-\widetilde{U}_{{F}}=g(e,Z)^{\top}\Delta_{\theta}\left(e_{\theta},Z\right).
	\end{aligned}
\end{equation}
Then it holds that
\begin{equation}\label{u1}
	\begin{aligned}
		\left|U_F\left(Z;\boldsymbol{\theta}\right)-\widetilde{U}_F\left(Z;\boldsymbol{\theta}\right)\right|
		&\leq \sup _{\theta \in \Theta} \mid \int_{\mathbf{Z}} \int_{e} g(e,Z)^{\top}\Delta\left(e_{\theta},Z\right) \times \left[d\hat{F}\left(e_{\theta}|Z\right)-d{F}\left(e_{\theta}|Z\right)\right]\cdot d\hat{G}(Z) \\
		&+\int_{\mathbf{Z}} \int_{e} g(e,Z)^{\top}\Delta\left(e_{\theta},Z\right)\cdot dF\left(e_{\theta}|Z\right)d\hat{G}(Z)\mid\\
		& \triangleq |B_1(h)|+|B_2(h)|.
	\end{aligned}
\end{equation}
where $\hat{F}$ and $\hat{G}$ are conditional and marginal empirical cumulative distribution functions respectively, and $F(u|Z)=\lim \limits_{n} \hat{F}(u|Z)$.

First, for $B_1(h)$,
since 
$$
\begin{aligned}
	g(e,Z)&\leq e\left[\left(\Phi\left(\frac{e+\alpha^{-1}}{h}\right)-I_{\{e> -\alpha^{-1}\}}\right)-\left(\Phi\left(\frac{e-\alpha^{-1}}{h}\right)-I_{\{e> \alpha^{-1}\}}\right)\right]\\
	&+\alpha^{-1}\left[\left(\Phi\left(\frac{e-\alpha^{-1}}{h}\right)-I_{\{e> \alpha^{-1}\}}\right)-\left(\Phi\left(\frac{e+\alpha^{-1}}{h}\right)-I_{\{e<-\alpha^{-1}\}}\right)+1\right]\\
	&\triangleq e\left[I_{11}+I_{12}\right]+\alpha^{-1}\left[I_{12}-I_{13}+1\right],
\end{aligned}
$$ 
where $I_{11}=\Phi\left(\frac{e+\alpha^{-1}}{h}\right)-I_{\{e> -\alpha^{-1}\}}$, $I_{12}=\Phi\left(\frac{e-\alpha^{-1}}{h}\right)-I_{\{e> \alpha^{-1}\}}$, and \\
$I_{13}=\Phi\left(\frac{e+\alpha^{-1}}{h}\right)-I_{\{e<-\alpha^{-1}\}}$.

Let $t=(e+\alpha^{-1})/h$, 
$g(e,Z)$ can be transformed into
$$
\begin{aligned}
	g(e,Z)\leq h\left\{\frac{e+\alpha^{-1}}{h}I_{11}-\frac{\alpha^{-1}}{h}I_{11}-\frac{e-\alpha^{-1}}{h}I_{12}-\frac{\alpha^{-1}}{h}I_{12}+\frac{\alpha^{-1}}{h}I_{12}-\frac{\alpha^{-1}}{h}I_{13}+\frac{\alpha^{-1}}{h}\right\}.
\end{aligned}
$$ 

Notice that $\sup _{t \in \mathbb{R}}|t(\Phi(t)-I(t>0))|=\sup _{t \in \mathbb{R}}\{t \Phi(-|t|) \operatorname{sgn}(t)\}$,
and it holds that $\lim _{t \rightarrow \infty}\{t \Phi(-|t|) \operatorname{sgn}(t)\}=0$. Under Condition C11, we conclude that
$$
\begin{aligned}
	g(e,Z)&\leq -\alpha^{-1}\left(I_{11}+I_{13}-1\right)+o(1) \\
	&=-2\alpha^{-1}\left[\Phi\left(\frac{e+\alpha^{-1}}{h}\right)-1\right].
\end{aligned}
$$
Moreover, it holds that 
$$
|B_{1}(h)|\leq \left|\int_{\mathbf{Z}} \int_{e} \Delta_{\theta}\alpha^{-1}2\left[\Phi\left(\frac{e+\alpha^{-1}}{h}\right)-1\right] \times \left[d\hat{F}\left(e_{\theta}|Z\right)-d{F}\left(e_{\theta}|Z\right)\right] d\hat{G}(Z)\right|.
$$
Let $\eta=\frac{r}{h}$, and integration by parts gives
$$
|B_{1}(h)|\leq2\alpha^{-1}\int_{\mathbf{Z}} \int_{e} \|\bar{\bar{\Delta}}_{\theta}\|\phi\left(\eta+\frac{\alpha^{-1}}{h}\right)\left[d\hat{F}\left(\eta h|Z\right)-d{F}\left(\eta h|Z\right)\right] d \eta d\hat{G}(Z).
$$
Under Condition C11, using the results on oscillations of empirical processes \citep{shorack86}, it holds that
$|B_{1}(h)|\leq O_p\left(\sqrt{h\log n\log\left(\frac{1}{h\log n}\right)}\right)$.

For $B_2(h)$, it holds that 
$$
|B_2(h)|\leq \int_{\mathbf{Z}} \int_{e} \|\bar{\bar{\Delta}}_{\theta}\|\alpha^{-1} 2 \left|\Phi\left(\frac{e+\alpha^{-1}}{h}\right)-1\right|
dF(e|Z)d\hat{G}(Z).
$$
Since 
\begin{equation}\label{tay}
	\begin{aligned}
		&\int_{e}\left[\Phi\left(\frac{e+\alpha^{-1}}{h}\right)-1\right]dF(e|Z)\\
		&=-\int_{\eta}F(\eta h|Z)\phi\left(\eta+\frac{\alpha^{-1}}{h}\right) d\eta
	\end{aligned}
\end{equation}
Apply second-order Taylor expansion to $F(\cdot|Z)$ at $h=0$, then 
$$
(\ref{tay})=-\frac{h^2}{2}\int_{\zeta}\zeta^2f^{\prime}(\zeta h^{*}|Z)\phi\left(\eta+\frac{\alpha^{-1}}{h}\right)d\eta,
$$
where $h^{*}\in (0,h)$, and $f^{\prime}(u|Z)=\partial^2 F(u|Z)/\partial u^2.$
Therefore, it holds that $|B_2(h)|\leq O_p(\alpha^{-1}h^2)=O_p(h^3\log n)$.

Combining the preceding arguments gives
$$
|U_F-\widetilde{U}_F|=O_p\left(\sqrt{h\log n\log\left(\frac{1}{h\log n}\right)}+h^3\right),
$$
under condition $nh^3\rightarrow 0$, it follows that $\tilde{U}_{F}(Z;\boldsymbol{\theta})-U_{F}(Z;\boldsymbol{\theta})=o_{p}(1)$ uniformly in $\boldsymbol{\theta}$.

{\bf Proof of Theorem 3 (3) and (4).} Denote $S_{p_m}=E_{F}\left[\dot{\widetilde{U}}_F\left(Z ; \widetilde{T}\left(F ; p_{m}\right)\right)\right]+\nabla^{2} p_{m}\left(|\boldsymbol{\theta}|;\lambda\right)$ for 
sequences $\left\{p_{m}\right\}_{m \geq 1}$ in $C^{\infty}(\Theta)$ converging to $p_n$ in $W^{2,2}(\Theta)$.
Thus
$$
\widetilde{\operatorname{IF}}_{p_{m}}(Z ; \widetilde T, F)=-S_{p_{m}}^{-1}\left({\widetilde{U}}_F\left(Z ; \widetilde T\left(F ; p_{m}\right)\right)+\nabla p_{m}\left(\widetilde{T}\left(F ; p_{m}\right)\right)\right),
$$
where $\dot{\widetilde{U}}$ is the derivative of $\widetilde{U}$.
Let 
For $p(\boldsymbol{\theta};\lambda)=\lambda|\boldsymbol{\theta}|$, 
there exists 
$$
g_{m}(t)=\frac{2}{m} \log \left(e^{t m}+1\right)-t \underset{m \rightarrow \infty}{\longrightarrow}|t|, 
$$
with 
$$
g_{m}^{\prime}(t):=\frac{2 e^{t m}}{e^{t m}+1}-1 \underset{m \rightarrow \infty}{\longrightarrow} \operatorname{sgn}(t),
$$
and
$$
g_{m}^{\prime \prime}(t)=\frac{2 m e^{t m}}{\left(e^{t m}+1\right)^{2}} \underset{m \rightarrow \infty}{\longrightarrow}\left\{\begin{array}{ll}
	0, & \text { if } t \neq 0 \\
	+\infty, & \text { otherwise },
\end{array}\right.
$$
which implies 
$\nabla^{2} p_{m}(T(F))_{j j} \rightarrow 0$ for $j=1,\ldots,s_n$ and $\nabla^{2} p_{m}(T(F))_{j j} \rightarrow \infty$ for $j=s_n+1,\ldots,p_n$.
Let 
$$
\dot{\widetilde{U}}_F=\left[\begin{array}{cc}
	\dot{{\widetilde{U}}}_{F_{11}} & \dot{{\widetilde{U}}}_{F_{12}} \\
	\dot{{\widetilde{U}}}_{F_{21}} & \dot{{\widetilde{U}}}_{F_{22}}
\end{array}\right]
$$
with $\dot{{\widetilde{U}}}_{F_{11}}$ being the $s_n \times s_n$ nonzero partitioned matrix. 
Then it holds that
$$
\begin{aligned}
	S_{p_{m}}^{-1}&=\left[\begin{array}{cc}
		E_{F}\left[\dot{{\widetilde{U}}}_{F_{11}}\right]+\nabla^{2} p_{m}^{1}(T(F)) & E_{F}\left[\dot{{\widetilde{U}}}_{F_{12}}\right] \\
		E_{F}\left[\dot{{\widetilde{U}}}_{F_{21}}\right] & E_{F}\left[\dot{{\widetilde{U}}}_{F_{22}}\right]+\nabla^{2} p_{m}^{2}(T(F))
	\end{array}\right]^{-1}\\
	&\underset{m \rightarrow \infty}{\longrightarrow}\left[\begin{array}{cc}
		\left(\widetilde{M}_{11}+P_{\lambda}\right)^{-1} & 0 \\
		0 & 0
	\end{array}\right]\triangleq {\Sigma}^{-1},
\end{aligned}
$$
where $\widetilde{M}_{11}=E_{F}\left[\dot{{\widetilde{U}}}_{11}(Z, \widetilde{T}(F))\right]$, $P_{\lambda}=\rm blockdiag\{p_{\lambda, j}^{\prime \prime}\left(\left|\theta_{j}^{*}\right|\right)\}$
for $j=1,\ldots,s_n$, and 0 otherwise. Thus we have 
\begin{equation}
	\begin{aligned}
		\widetilde{\operatorname{IF}}(Z ; \widetilde T, F)&=\lim_{m\rightarrow \infty} \widetilde{\operatorname{IF}}_{p_{m}}(Z ; \widetilde T, F)\\
		&={\Sigma}^{-1}\left(\widetilde{U}_F\left(Z;\widetilde{T}(F)\right)+\nabla p\left(\widetilde{T}\left(F ; p\right)\right)\right).
	\end{aligned}
\end{equation}

Using Lemma, $|U_{F_{11}}-\widetilde{U}_{F_{11}}|\rightarrow 0$ a.s., thus
$$
\widetilde{\operatorname{IF}}(Z ; F, \widetilde{T}) \xrightarrow[n \rightarrow \infty]{P} \mathrm{IF}(Z ; F, T)=-S^{-1}\left(U_{F}+p^{\prime}(|\boldsymbol{\theta}|)\right),
$$
where
$S^{-1}=\rm blockdiag \{(M_{11}+P_{\lambda})^{-1},0\}$ with $M_{11}=E_{F}\left[\dot{U}_{F_{11}}(Z, T(F))\right]$.
Thus the claimed results (3) and (4) can be proved.

It is noteworthy that for sufficiently small $\alpha^{-1}$, the smoothed estimating equation $\widetilde{U}_F$ exhibits boundedness concerning the outliers present in either the response or covariate domains, irrespective of the bandwidth function $h(n)$. This indicates that our proposed estimators possess strong robustness against both the outliers in response and outlying influence points.

\section{The convergence analysis on Bi-convex optimization (2.2)}

Recall the biconvex optimization problem (2.2), 
\begin{equation}\label{obj}
	(\widehat{\boldsymbol{\beta}}, \widehat{\mathbf{w}})=\underset{\boldsymbol{\beta}, \mathbf{w}}{\arg \min }f(\boldsymbol{\beta}, \mathbf{w}):=\left\{\frac{1}{2} \sum_{i=1}^{n} 
	\ell_{\alpha}\left(w_{i}(y_{i}-x_{i}^{\top} \boldsymbol{\beta})\right) +\mu \sum_{i=1}^{n} {\varpi}_{i}\left|1-w_{i}\right|+\lambda \sum_{j=1}^{p_n} \left|\beta_{j}\right|\right\},
\end{equation}
over the joint parameter space $(\boldsymbol{\beta},\mathbf{w})\in(\Theta_1=\mathbf{R}^{p_n},\Theta_2=(0,1)^{n})$.
The Alternative convex search algorithm is employed for solving it\citep{Gor07}, by splitting (2.2) into inner layer(for updating $\widehat{\boldsymbol{\beta}}$) and outer layer(for updating $\widehat{\mathbf{w}}$).
The final solution path $(\widehat{\boldsymbol{\beta}}, \widehat{\mathbf{w}})$ is a partial optima, satisfying: 
\begin{equation}
	f(\widehat{\boldsymbol{\beta}}, \widehat{\mathbf{w}})\leq f({\boldsymbol{\beta}}, \widehat{\mathbf{w}}),
	f(\widehat{\boldsymbol{\beta}}, \widehat{\mathbf{w}})\leq f(\widehat{\boldsymbol{\beta}}, {\mathbf{w}})
\end{equation}
for any $\boldsymbol{\beta} \in \Theta_1$, and $\mathbf{w} \in \Theta_2$.

For the sequence joint estimator $z_k=\left(\widehat{\boldsymbol{\beta}}^{(k)},\widehat{\mathbf{w}}^{(k)}\right)$, convergence to a unique estimator is not guaranteed.
But the sequence does have at least one accumulation point, which may not necessarily be a local optimum. 
Nevertheless, Theorem 4.9 in \citet{Gor07} shows it is at least a optima.
However, when the optima $\widehat{\mathbf{w}}$ exists and is unique, and $f(\boldsymbol{\beta}, \mathbf{w})$ is strictly convex w.r.t $\boldsymbol{\beta}$ when given $\mathbf{w}$, and vice versa. The sequence $\left\{z_k\right\}$ will converge to a unique point.

\subsection{Convergence analysis on the outer layer iteration}\label{w}

In the outer layer, we apply greedy coordinate descent to solve convex optimization:
\begin{equation}\label{eqw1}
	\widehat{\mathbf{w}}=\arg \min _{\mathbf{w}}\left\{\sum_{i=1}^{n}\ell_{\alpha}\left(w_{i}(y_{i}-x_{i}^{\top} \boldsymbol{\beta})\right)+\mu \sum_{i=1}^{n} {\varpi}_{i}\left|1-w_{i}\right|\right\},
\end{equation}
and let $\tilde{w}_i=1-w_i$, then (\ref{eqw1}) can be transformed into 
\begin{equation}\label{eqw2}
	\hat{\tilde{w}}_i=\arg\min_{\tilde{w}_i} l(\tilde{w}_i):= \ell_{\alpha}\left((1-\tilde{w}_i)r_i\right)+\bar{\mu_i}|\tilde{w}_i|,	
\end{equation}
where $\bar{\mu}_i=\mu{\varpi}_{i}$, and 
$$
\ell_{\alpha}(\left((1-\tilde{w}_i)r_i\right))=\left\{
\begin{array}{lll}
	2 \alpha^{-1}|e_i|-\alpha^{-2} & \text {if } & |e_i|>\alpha^{-1} \\
	\left[(1-\tilde{w}_i)r_i\right]^{2} & \text { if } & |e_i| \leq \alpha^{-1}
\end{array} \right.
$$
with $r_i=y_{i}-x_{i}^{\top} \boldsymbol{\beta}$ and $e_i=r_i-r_i\tilde{w}_i$.
Differentiating (\ref{eqw2}) w.r.t $\tilde{w}_i$ gives:
$$
\begin{aligned}
	\partial{l(\tilde{w}_i)}/\partial{\tilde{w}_i}&=\bar{\mu}_i-2\alpha^{-1}r_i \sgn\left(r_i(1-\tilde{w}_i)\right)\\
	&\leq \bar{\mu_i}-\ell_{\alpha}\left(r_i\right).
\end{aligned}
$$
Thus, when $\ell_{\alpha}\left(r_i\right)>\bar{\mu}_i$, the gradient is less than 0, hence we update $\tilde{w}_i$, otherwise, $\tilde{w}_i\leftarrow 1$.

The aforementioned coordinate descent algorithm is weakly consistent when following regularity conditions hold:

S1. The parameter vector $\mathbf{w}$ is confined to a compact domain $K \subset (0,1)^n$.
The true parameter vector $\mathbf{w}^{*}$ is an interior point of $K$.

S2. In heterogeneous model, $y_{i}=x_{i}^{\top} \boldsymbol{\beta}+\varepsilon_{i}$, $\varepsilon_{i}=g\left(x_{i}\right) \tilde{\varepsilon}_{i}$ with 
$\tilde{\varepsilon}_{i}$ being independent of $x_i$. The weighted random errors $e_i=w_i(y_i-x_{i}^{\top} \boldsymbol{\beta})$ are independent, and the distribution function (denoted as $F_i(\cdot)$) of $\frac{e_i}{g(x_i)}$ is symmetric at $0$, i.e. $F_i(0)=1/2$.

S3. For every $c>0$ there exists an $f>0$ with 
$$
\inf_i \min \left\{F_i(c)-\frac{1}{2}, \frac{1}{2}-F_i(-c)\right\} \geq f.
$$

S4. The residual $r_i$ satisfy $\left\|r_i\right\|_2 \leq B$ for some $B\geq 0$.

S5. For some $t>0$ and $d>0$, the $r_i$ satisfy$\inf\frac{1}{n} \sum_{i=1}^n 1_{\left\{\left|r_i\right| \geq t\right\}} \geq d,$
for sufficiently large $n$.

\text{\bf Remark.} It is straightforward to verify Conditions S1 and S4 over the parameter space $\Theta_2$ in equation (\ref{eqw2}) provided that Condition C4 in section 3.2 is satisfied. The development of Condition S2 builds upon and extends the work of \citet{gao16}, who considered the heterogeneous setting but in the context of variance function linear models. Similarly, Conditions S3 and S5 are inspired by the work of \citet{wu08}, and our subsequent proof follows the same spirit as their work. Nonetheless, our approach generalizes their results to the case of heterogeneous random errors.

{\bf Theorem 4.} Suppose Conditions S1--S5 hold, the sequence of estimators $\widehat{\tilde{\mathbf{w}}}$ minimizing (\ref{eqw2})
given $\boldsymbol{\beta}$, is weakly consistent. As a result, the convergence of sequence $\widehat{\mathbf{w}}$ minimizing (\ref{eqw1}) can be derived directly.

{\bf Proof.} For the sake of parsimony, we define parameter vector $\theta:=\tilde{\mathbf{w}}=1-\mathbf{w}$ with $i$th component $\theta_i$.
Denote
\begin{equation}\label{eqw3}
	f_n(\theta)=g_n(\theta)+\sum_{j=1}^n\bar{\mu}_j|\theta_j|,
\end{equation}
where $g_n(\theta)=\sum_{i=1}^n\ell_{\alpha}\left(r_i-r_i\theta_i\right)$. Denote $\theta^{*}$ as the true parameter vector, and an optima of (\ref{eqw3})
is $\widehat{{\theta}}$. Let $d_n(\theta)=f_n(\theta)-f_n(\theta^{*})$, then it holds that $d_n(\widehat{\theta})\leq 0$.

Our goal is to prove $\widehat{{\theta}}$ coverges to $\theta^{*}$ a.s. The proof is completed based on three claims:

\Rmnum{1}. Denote compact set $C\subset K=(0,1)^n$, for any estimate $\theta_i \in(0,1)$, it can be directly derived
$$
\left|d_n(\theta)-d_n\left(\widehat{\theta}\right)\right| \leq R\left\|\theta-\widehat{\theta}\right\|_2+ \sum_{j=1}^n \bar{\mu}_i\left|\theta_j-\widehat{\theta}_j\right|,
$$
where $R=2\alpha^{-1}B$. Hence, the uniform continuity of $d_n(\theta)$ holds for all $\widehat{\theta}$ in a neighborhood of $\theta$.

\Rmnum{2}. $\inf _{\theta \in C} \mathrm{E}\left[d_n(\theta)\right]\geq \eta$, for $\eta>0$ and large $n$.

This fact gives
$$
\begin{aligned}
	\operatorname{Pr}\left[d_n(\theta) \geq \eta-\delta\right]& \geq \operatorname{Pr}\left[d_n(\theta)\geq \inf_{\theta\in C}\mathrm{E}\left[d_n(\theta)\right]-\delta\right]\\
	&\geq 1-\frac{\var\left[d_n(\theta)\right]}{\delta^2}
\end{aligned}
$$
The second inequality holds using Chebshev inequality. Next, let $\delta=\frac{\eta}{2}$. 
$$
\begin{aligned}
	\operatorname{Pr}\left[\inf _{\widehat{\theta} \in \mathcal{B}(\theta)} d_n\left(\widehat{\theta}\right) \geq \frac{1}{4} \eta\right] &\geq \operatorname{Pr}\left[\left\{\inf \left(d_n\left(\widehat{\theta}\right)-d_n(\theta)\right) \geq-\frac{1}{4} \eta\right\}\bigcup\left\{d_n(\theta)\geq \frac{\eta}{2}\right\}\right]\\
	&= \operatorname{Pr}\left(d_n(\theta)\geq \frac{\eta}{2}\right)\\
	& \geq 1-\frac{\var\left[d_n(\theta)\right]}{(\eta/2)^2}.
\end{aligned}
$$
The first equality holds by using fact \Rmnum{1}. By using Chebshev inequality again, we can derive the second inequality.
When we take $C=\mathcal{B}^{c}({\theta}^*)\subset K$, the complement of the neighborhood of $\theta^*$,
according to Borel-Lebesgue covering theorem, there is a set of finite cover set $N_i$ for $i=1,\ldots,m$ s.t. $C\subset  \bigcup_{i=1}^{m} N_i$, then it holds that
$$
\begin{aligned}
	\operatorname{Pr}\left[\inf _{\widehat{\theta} \in C} d_n\left(\widehat{\theta}\right)<\frac{1}{4} \eta\right]& \leq \sum_{i=1}^m \operatorname{Pr}\left[\inf _{\widehat{\theta} \in N_i} d_n\left(\widehat{\theta}\right)<\frac{1}{4} \eta\right]\\
	& \leq \sum_{i=1}^m \frac{\var \left[d_n(\theta)\right]}{(\eta/2)^2}
\end{aligned}
$$

\Rmnum{3}.  $\lim _{n \rightarrow \infty} \sup _{\theta \in C} \var\left[d_n(\theta)\right]=0$.

By using this, $\frac{\var\left[d_n(\theta)\right]}{(\eta/2)^2}<\varepsilon$ for any $\varepsilon>0$.
Then for any $\widehat{\theta}\in C$, $d_n(\widehat{\theta})=f_n(\widehat{\theta})-f_n({\theta}^*)\rightarrow 0$. 
This implies for any estimate $\widehat{{\theta}} \in K$, the probability of $\widehat{{\theta}}\in C$ is at most $\varepsilon$, whereas
with high probability, $\widehat{{\theta}}\in \mathcal{B}(\theta^{*})$, where $d_n(\widehat{{\theta}})\rightarrow 0$, thus any optima of (\ref{eqw3})
$\widehat{{\theta}}$ converges to the true value $\theta^{*}$ a.s. 

The completion of the proofs for\Rmnum{2}-\Rmnum{3} can be accomplished by following similar steps as in \citet{wu08}, with the caveat that modified versions of Conditions S2--S5 are needed.

\subsection{Convergence analysis on the inner layer iteration}
Given estimate $\mathbf{w}$, that has been discussed in section \ref{w}, the optima of (2.4), denoted as $\widehat{\boldsymbol{\beta}}$, can be derived from
\begin{equation}
	\widehat{\boldsymbol{\beta}}=\arg \min _{\boldsymbol{\beta}}\left\{\sum_{i=1}^n \ell_\alpha(w_i\left(y_i-x_i^{\top}\boldsymbol{\beta}\right))+\lambda \sum_{j=1}^{p_n}\left|\beta_j\right|\right\},
\end{equation}
and we denote $L_n(\boldsymbol{\beta})=\sum_{i=1}^n \ell_\alpha(w_i\left(y_i-x_i^{\top}\boldsymbol{\beta}\right))$, then the updated estimate at $k$th step is:
\begin{equation}
	\widehat{\boldsymbol{\beta}}^{(k)}=\underset{\|\boldsymbol{\beta}\|_{1} \leq \rho}{\operatorname{argmin}}\left\{L_n\left(\widehat{\boldsymbol{\beta}}^{(k-1)}\right)+\left[\nabla L_n\left(\widehat{\boldsymbol{\beta}}^{(k-1)}\right)\right]^T\left(\boldsymbol{\beta}-\widehat{\boldsymbol{\beta}}^{(k-1)}\right)+\frac{\gamma_n}{2}\left\|\boldsymbol{\beta}-\widehat{\boldsymbol{\beta}}^{(k-1)}\right\|_2^2+\lambda\|\boldsymbol{\beta}\|_1\right\}
\end{equation}
by using local isotropic quadratic approximation, and parameter $\gamma_n$ is related to second-order smoothness of $L_n(\boldsymbol{\beta})$, and is independent of 
$\alpha$. The convergence rate
$$	
\begin{aligned}
	\left\|\widehat{\boldsymbol{\beta}}^{(k)}-\boldsymbol{\beta}^*\right\|_2 & \leq\left\|\widehat{\boldsymbol{\beta}}^{(k)}-\widehat{\boldsymbol{\beta}}\right\|_2+\left\|\widehat{\boldsymbol{\beta}}-\boldsymbol{\beta}^*\right\|_2.
\end{aligned}	
$$
By using Lemma 4 and Theorem 4 in \citet{fan17}, the first term converges to 0 a.s. with large enough sample size. The convergence rate of the second term 
is 
$$
O\left(\alpha_0^{-1}\left(d_n+q_n\right)^{1 / 2}\left(\frac{\log \left(p_n+n\right)}{n}\right)^{\min \{\delta /(1+\delta), 1 / 2\}}\right),
$$ as shown in (3.10)
in Theorem 2. Hence, the convergence rate of algorithm (2.4) is in the same order as that of second term.

\section{Finite-sample breakdown point}

Since $L_1$ norm $\|\theta\|_1$ and the Euclidean norm $\|\theta\|_2$ are topologically equivalent, there exists a constant $c_1>0$ and $c_2>0$ such that
$c_1\|\theta\|_2 \leq \|\theta\|_1 \leq c_2 \|\theta\|_2$. Denote
$$
Q(\theta)=\sum_{i=1}^n \ell_{\alpha}\left(y_i-z_{i,\beta}^{\top}\theta\right)+n \lambda \|\theta\|_1.
$$
Suppose we contaminate $m=0$ observation, namely all sample are good observations.
First, we prove that the breakdown point is at least 1/n.

For an initial estimate $\widetilde{\theta}$ such that $\|\widetilde{\theta}\|_2\leq M$, $M_y=\max \limits_{i} |y_i|$, and $M_z=\max \limits_{i} \|z_{i,\widetilde{\beta}}\|_{\infty}$.
It holds that
$$
\begin{aligned}
	Q\left(\widetilde{\theta}\right)&=\sum_{i=1}^n \ell_{\alpha}\left(y_i-z_{i,\widetilde{\beta}}^{\top}\widetilde{\theta}\right)+\lambda n \|\widetilde{\theta}\|_1\\
	&\leq n \alpha^{-2} +n\left[2\alpha^{-1}M_y-\alpha^{-2}\right]+\lambda n c_2 \|\widetilde \theta\|_2 \\
	&=2n\alpha^{-1}M_y+\lambda n c_2 \|\widetilde{\theta}\|_2 \\
	& \leq 2n\alpha^{-1} M_y + \lambda n c_2 M \\
	&\triangleq M^{\star}.
\end{aligned}
$$
For any vector $\theta \in \mathbb{R}^{p_n+n}$, when $\|\theta\|_2 \geq \frac{M^{\star}}{\lambda c_1}$, it holds that
$$
Q\left(\theta\right) \geq \lambda \|\theta\|_1 \geq \lambda c_1 \|\theta\|_2=M^{\star} \geq Q\left(\widetilde{\theta}\right).
$$
Since $Q\left(\widehat{\theta}\right)\leq Q\left(\widetilde{\theta}\right)$, it indicates the breakdown point is at least $\frac{1}{n}$. 

Second, we conclude that $\mathrm{BP}\left(\widehat{\boldsymbol{\theta}} ;\widetilde{\mathbf{Z}}^{m}_{(n)} \right)\leq \frac{1}{n}$.
Motivated by the proof of Theorem 1 in \citet{Alf13}, for the parameter vector 
$\boldsymbol{\theta}_{\gamma}=(\gamma, 0, \ldots, 0)^{\prime} \in \mathbb{R}^{p_n+n}$. 
We contaminate $m=1$ observation, denote the contaminated sample 
$\mathbf{Z}_{\gamma,\tau}=\left(z(\tau)^{\prime}, y(\gamma, \tau)\right)^{\prime}=((\tau, 0, \ldots, 0), \gamma \tau)^{\prime}$, for some $\tau, \gamma>0$.
Assume that there exists a constant $M$ such that
$\sup \limits_{\tau, \gamma}\left\|\hat{\theta}\left(\mathbf{Z}_{\gamma, \tau}\right)\right\|_{2} \leq M$.
Let $\gamma=M+2$. 
Denote $M_{z_1}=\max\limits_{j} \|z_{1j,\beta}\|_{\infty}$.
Define $\tau>0$ such that 
$$
\begin{aligned}
	Q\left(\theta_{\gamma}\right)&=\sum_{i=1}^n \ell_{\alpha}\left(y_i-z_{i,\beta}^{\top}\theta_{\gamma}\right)+n\lambda \|\theta\|_1\\
	& \leq (n-1)\left[2\alpha^{-1}(M_y+\gamma M_{z_1})-\alpha^{-2}\right]+n\lambda \gamma \\
	& \leq (n-1)\left[2\alpha^{-1}(M_y+\gamma M_{z_1})\right]+n \lambda \gamma \\
	& \leq \ell_{\alpha}(\tau).
\end{aligned}
$$
For $\|\theta\|_2\leq \gamma-1$, it holds that $\theta_1 \leq \|\theta\|_2\leq \gamma -1$.
Therefore,
$$
\begin{aligned}
	Q\left(\theta\right)&=\sum_{i=1}^n \ell_{\alpha}\left(y_i-z_{i,\beta}^{\top}\theta\right)+n\lambda \|\theta\|_1\\
	& \geq \ell_{\alpha}\left(y(\gamma,\tau)-z(\tau)^{\top}\theta\right)\\
	& \geq \ell_{\alpha} \left(\gamma \tau-\tau \theta_1\right)\\
	& \geq \ell_{\alpha}(\tau)\\
	& \geq Q\left(\theta_{\gamma}\right).
\end{aligned}
$$
Since $Q\left(\widehat{\theta}\right)\leq Q\left(\theta_{\gamma}\right)$, thus $\|\widehat{\theta}\|_2\geq \gamma-1$, 
also $\sup \limits_{\gamma} \|\widehat{\theta}\left(Z_{\gamma},\gamma\right)\|_2 \geq \gamma-1=M+1$, which indicates that 
the breakdown point is no more than 1/n. To conclude, we deduce the breakdown point of estimator $\theta$.

\newpage

\section{Additional Results}

\subsection{Simulation studies}

We repeat simulations \Rmnum{1}-\Rmnum{2} in Section 4.1, adhering to the previously established settings, except here we take $c = 0.2, 0.3$. The corresponding results are presented in Tables 1--4.

\begin{footnotesize}
	\begin{longtable}{lcccccccccc}
		\caption{Comparison of PM\citep{kong18}, PWLAD, AHuber, PIQ and the proposed method PWHL
			under three cases in homogeneous model for $c=0.2$.}\label{suptab1-2}\\
		\hline
		(T, Case)  & Method & M & S & JD &  EE & EE$_{non}$  & FZR & FPR & SR & CR \\
		\hline
		\endfirsthead	
		\hline
		(T, Case)  & Method & M & S & JD &  EE & EE$_{non}$  & FZR & FPR & SR & CR \\
		\hline
		\endhead
		\hline
		\endfoot
		(T1, Case 1) & PM & 0.2010 & \textbf{0.0092} & 0.02 & \textbf{0.5440} & 11.5742 & 0.0333 & 0.0011 & 0.62 & 0.91 \\
		& AHuber & - & - & -& 1.0354 & 6.7226 & 0.4333 & 0.0020 & 0.13 & 0.31 \\ 
		& PWLAD & 0.1050 & 0.2988 & 0.10 & 0.6625 & 12.8982 & \textbf{0.0000} & 0.0715 & 0.00 & \textbf{1.00} \\  
		& PWHL & \textbf{0.0950} & 0.1812 & \textbf{0.10} & {0.5983} & 0.5983 & 0.0333 & \textbf{0.0000} & \textbf{0.90} & 0.90 \\ 
		& PIQ & 0.1300 & 0.1575 & 0.00 & 1.0619 & 0.8956 & 0.4000 & 0.0055 & 0.00 & 0.10 \\
		\hline
		(T1, Case 2) & PM & 0.9980 & 0.0001 & 0.00 & \textbf{0.2944} & 14.3258 & 0.0000 & 0.0112 & 0.03 & 1.00 \\ 
		& AHuber & - & - & -  & 1.3587 & 1.4439 & 0.0000 & 0.2480 & 0.00 & 1.00 \\ 
		& PWLAD & 0.0200 & 0.1812 & 0.80 & 0.5958 & 13.6234 & 0.0000 & 0.0756 & 0.00 & 1.00 \\ 
		& PWHL & \textbf{0.0000} & \textbf{0.0000} & \textbf{1.00} & 0.4030 & \textbf{0.4009} & 0.0000 & \textbf{0.0003} & \textbf{0.90} & 1.00 \\
		& PIQ & 0.0000 & 0.1250 & 1.00 & 1.0376 & 0.8864 & 0.3667 & 0.0053 & 0.00 & 0.20 \\
		\hline
		(T1, Case 3) & PM & 0.0000 & 0.0002 & 1.00 & \textbf{0.3823} & 13.1415 & 0.0000 & 0.0007 & 0.74 & 1.00 \\
		& AHuber & - & - & -& 1.3856 & 1.3856 & 1.0000 & 0.0000 & 0.00 & 0.00 \\ 
		& PWLAD & 0.0000 & 0.2238 & 1.00 & 0.5946 & 14.8172 & 0.0000 & 0.0766 & 0.00 & 1.00 \\ 
		& PWHL & 0.0000 & \textbf{0.0000} & 1.00 & 0.3992 & \textbf{0.3992} & \textbf{0.0000} & \textbf{0.0000} & \textbf{1.00} & 1.00 \\
		& PIQ & 0.0000 & 0.1250 & 1.00  & 7.1137 & 1.3856 & 1.0000 & 0.0101 & 0.00 & 0.00 \\
		\hline
		(T2, Case 1) & PM & 0.3090 & \textbf{0.0338} & 0.00 & 0.8385 & 9.3089 & 0.3033 & \textbf{0.0008} & 0.33 & 0.53 \\
		& AHuber & - & - & - & 1.1186 & 6.4349 & 0.4967 & 0.0043 & 0.04 & 0.27 \\ 
		& PWLAD & 0.1250 & 0.3025 & \textbf{0.20} & 0.9765 & 12.9784 & \textbf{0.0667} & 0.0723 & 0.00 & 0.80 \\
		& PWHL & \textbf{0.1050} & 0.2775 & 0.10 & \textbf{0.8038} & \textbf{0.8006} & 0.1000 & 0.0018 & \textbf{0.40} & \textbf{0.80} \\
		& PIQ & 0.1525 & 0.1631 & 0.00 & 1.2636 & 1.0461 & 0.5600 & 0.0068 & 0.00 & 0.07 \\
		\hline
		(T2, Case 2) & PM & 0.9765 & 0.0115 & 0.00 & \textbf{0.4377} & 13.7905 & 0.0067 & 0.0109 & 0.01 & 0.98 \\ 
		& AHuber & - & - & - & 1.3548 & 1.4677 & 0.0100 & 0.2505 & 0.00 & 0.99 \\ 
		& PWLAD & 0.0100 & 0.2812 & 0.80 & 0.6934 & 14.2653 & 0.0000 & 0.0640 & 0.00 & 1.00 \\
		& PWHL & 0.0000 & \textbf{0.0013} & 1.00 & \textbf{0.6448} & \textbf{0.6444} & \textbf{0.0000} & \textbf{0.0008} & \textbf{0.70} & 1.00 \\
		& PIQ & 0.0000 & 0.1250 & 1.00 & 1.0616 & 0.9040 & 0.4500 & 0.0059 & 0.00 & 0.21 \\
		\hline
		(T2, Case 3) & PM & 0.0000 & 0.0251 & 1.00 & \textbf{0.5365} & 11.9988 & 0.0400 & 0.0008 & 0.64 & 0.89 \\
		& AHuber & - & - & - & 1.3856 & 1.3856 & 1.0000 & 0.0000 & 0.00 & 0.00 \\
		& PWLAD  & 0.0000 & 0.2575 & 1.00 & 0.7635 & 13.7010 & 0.0000 & 0.0675 & 0.00 & 1.00 \\ 
		& PWHL & 0.0000 & \textbf{0.0025} & 1.00 & 0.5645 & \textbf{0.5635} & 0.0000 & 0.0013 & \textbf{0.80} & 1.00 \\
		& PIQ & 0.0000 & 0.1250 & 1.00 & 7.1313 & 1.4928 & 0.9833 & 0.0099 & 0.00 & 0.00 \\
		\hline
	\end{longtable}	
\end{footnotesize}

\begin{footnotesize}
	\begin{longtable}{lcccccccccc}
		\caption{Comparison of PM\citep{kong18}, PWLAD, AHuber, PIQ and the proposed method PWHL
			under three cases in homogeneous model for $c=0.3$.}\label{tab1-3}\\
		\hline
		(T, Case)  & Method & M & S & JD &  EE & EE$_{non}$  & FZR & FPR & SR & CR \\
		\hline
		\endfirsthead	
		\hline
		(T, Case)  & Method & M & S & JD &  EE & EE$_{non}$  & FZR & FPR & SR & CR \\
		\hline
		\endhead
		\hline
		\endfoot
		(T1, Case 1) & PM & 0.2080 & \textbf{0.0086} & 0.00 & 0.5901 & 10.8620 & 0.0533 & 0.0012 & 0.56 & 0.88 \\
		& AHuber  & - & - & - & 1.3087 & 2.9715 & 0.8500 & 0.0005 & 0.00 & 0.02 \\ 
		& PWLAD & 0.1200 & 0.3700 & \textbf{0.20} & 0.8247 & 12.2333 & \textbf{0.0000} & 0.0713 & 0.00 & \textbf{1.00} \\  
		& PWHL & \textbf{0.0733} & 0.1971 & 0.10 & \textbf{0.5280} & \textbf{0.5135} & 0.0333 & 0.0010 & \textbf{0.60} & 0.90 \\
		& PIQ & 0.1367 & 0.2586 & 0.00 & 1.5548 & 1.2878 & 0.8667 & 0.0091 & 0.00 & 0.00 \\
		\hline
		(T1, Case 2) & PM & 0.9977 & 0.0003 & 0.00 & \textbf{0.2987} & 14.3303 & 0.0000 & 0.0116 & 0.01 & 1.00 \\
		& AHuber & - & - & - & 1.3614 & 1.4356 & 0.0200 & 0.2397 & 0.00 & 0.98 \\ 
		& PWLAD & 0.0067 & 0.1843 & 0.80 & 0.5361 & 14.4119 & 0.0000 & 0.0610 & 0.00 & 1.00 \\
		& PWHL & 0.0067 & \textbf{0.0000} & 0.80 & 0.3559 & \textbf{0.3522} & \textbf{0.0000} & \textbf{0.0003} & \textbf{0.90} & 1.00 \\
		& PIQ & \textbf{0.0000} & 0.2000 & \textbf{1.00} & 1.1947 & 1.0528 & 0.5667 & 0.0068 & 0.00 & 0.00 \\
		\hline
		(T1, Case 3) & PM & 0.3510 & 0.0307 & 0.64 & 3.4056 & 16.0464 & 0.2000 & 0.0303 & 0.35 & 0.62 \\ 
		& AHuber & - & - & - & 1.3856 & 1.3856 & 1.0000 & 0.00 & 0.00 & 0.00 \\ 
		& PWLAD & 0.0000 & 0.1743 & 1.00 & 0.5324 & 14.1476 & \textbf{0.0000} & 0.0584 & 0.00 & \textbf{1.00} \\ 
		& PWHL & 0.0000 & \textbf{0.0000} & 1.00 & \textbf{0.4896} & \textbf{0.4896} & 0.0667 & \textbf{0.0000} & \textbf{0.80} & 0.80 \\
		& PIQ & 0.0000 & 0.2000 & 1.00 & 6.8924 & 1.4093 & 0.9900 & 0.0100 & 0.00 & 0.00 \\
		\hline
		(T2, Case 1) & PM & 0.3560 & \textbf{0.0216} & 0.00 & 0.9082 & 8.6033 & 0.3367 & 0.0007 & 0.34 & 0.51 \\
		& AHuber  & - & - & -& 1.3368 & 2.3972 & 0.9000 & \textbf{0.0004} & 0.01 & 0.02 \\ 
		& PWLAD & 0.1600 & 0.2814 & 0.10 & 1.1452 & 12.0748 & \textbf{0.0333} & 0.0761 & 0.00 & \textbf{0.90} \\ 
		& PWHL & 0.1333 & 0.1914 & 0.10 & \textbf{0.8607} & \textbf{0.8182} & 0.1667 & 0.0086 & \textbf{0.30} & 0.60 \\
		& PIQ & \textbf{0.1250} & 0.2536 & 0.00 & 1.4969 & 1.2286 & 0.7733 & 0.0084 & 0.00 & 0.00 \\
		\hline
		(T2, Case 2) & PM & 0.9833 & 0.0149 & 0.00 & \textbf{0.4311} & 13.6801 & \textbf{0.0033} & 0.0116 & 0.00 & 0.99 \\ 
		& AHuber & - & - & - & 1.3617 & 1.4346 & 0.0100 & 0.2426 & 0.00 & \textbf{0.99} \\ 
		& PWLAD & 0.0033 & 0.2000 & 0.90 & 0.7858 & 13.0660 & 0.0333 & 0.0640 & 0.00 & 0.90 \\
		& PWHL & 0.0000 & \textbf{0.0000} & 1.00 & 0.6732 & \textbf{0.6730} & 0.0667 & \textbf{0.0015} & \textbf{0.50} & 0.90 \\
		& PIQ & 0.0000 & 0.2000 & 1.00 & 1.3791 & 1.1595 & 0.6833 & 0.0077 & 0.00 & 0.01 \\
		\hline
		(T2, Case 3) & PM & 0.3290 & 0.0421 & 0.66 & 3.4261 & 13.9035 & 0.3567 & 0.0289 & 0.17 & 0.47 \\ 
		& AHuber & - & - & -  & 1.3856 & 1.3856 & 1.0000 & 0.0000 & 0.00 & 0.00 \\
		& PWLAD & 0.0000 & 0.2300 & 1.00 & 0.7636 & 13.5683 & 0.0000 & 0.0617 & 0.00 & 1.00 \\ 
		& PWHL & 0.0000 & \textbf{0.0014} & 1.00 & \textbf{0.6202} & \textbf{0.6154} & \textbf{0.0000} & \textbf{0.0023} & \textbf{0.30} & \textbf{1.00} \\
		& PIQ & 0.0000 & 0.2000 & 1.00 & 6.5919 & 1.5087 & 0.9900 & 0.0100 & 0.00 & 0.00 \\
		\hline
	\end{longtable}	
\end{footnotesize}	

\begin{footnotesize}
	\begin{longtable}{lcccccccccc}
		\caption{Comparison of PM\citep{kong18}, PWLAD, AHuber, PIQ and the proposed method PWHL
			under three cases in heterogeneous model for $c=0.2$.}\label{tab2-2}\\
		\hline
		(T, Case)  & Method & M & S & JD &  EE & EE$_{non}$  & FZR & FPR & SR & CR \\
		\hline
		\endfirsthead	
		\hline
		(T, Case)  & Method & M & S & JD &  EE & EE$_{non}$  & FZR & FPR & SR & CR \\
		\hline
		\endhead
		\hline
		\endfoot
		(T1, Case 1) & PM & 0.4180 & \textbf{0.1504} & 0.01 & 1.0959 & 5.9143 & 0.6200 & 0.0007 & 0.12 & 0.24 \\ 
		& AHuber & - & - &- & 1.9035 & 3.8049 & 0.8300 & 0.0213 & 0.00 & 0.04 \\ 
		& PWLAD  & 0.2000 & 0.3275 & 0.00  & 1.3036 & 12.2169 & 0.1000 & 0.0746 & 0.00 & 0.70 \\ 
		& PWHL & \textbf{0.1350} & 0.3100 & 0.00 & \textbf{0.5951} & \textbf{0.5943} & \textbf{0.0333} & \textbf{0.0003} & \textbf{0.80} & \textbf{0.90} \\
		& PIQ & 0.3000 & 0.2000 & 0.00 & 1.3560 & 1.0250 & 0.5000 & 0.0063 & 0.00 & 0.10 \\
		\hline
		(T1, Case 2) & PM & 0.7365 & 0.2534 & 0.00 & \textbf{0.4135} & 13.5798 & 0.0267 & 0.0112 & 0.02 & 0.95 \\ 
		& AHuber & - & - & - & 1.3556 & 1.4841 & 0.1200 & 0.2866 & 0.00 & 0.88 \\ 
		& PWLAD & 0.0050 & 0.3550 & 0.90 & 1.1004 & 12.3443 & 0.1000 & 0.0625 & 0.00 & 0.80 \\  
		& PWHL & 0.0050 & 0.2650 & 0.90 & 0.5372 & \textbf{0.4390} & \textbf{0.0000} & 0.0065 & \textbf{0.30} & \textbf{1.00} \\
		& PIQ & \textbf{0.0000} & \textbf{0.1250} & \textbf{1.00} & 1.6093 & 0.8996 & 0.4000 & \textbf{0.0055} & 0.00 & 0.20 \\
		\hline
		(T1, Case 3) & PM & 0.0000 & 0.1771 & 1.00 & 0.9934 & 7.8545 & 0.4800 & 0.0009 & 0.16 & 0.30 \\ 
		& AHuber & - & - & - & 1.6441 & 2.3245 & 0.9367 & 0.0144 & 0.00 & 0.06 \\  
		& PWLAD & 0.0000 & 0.3500 & 1.00 & 0.9945 & 10.4177 & \textbf{0.1000} & 0.0602 & 0.00 & 0.70 \\
		& PWHL & 0.0000 & 0.1375 & 1.00 & \textbf{0.7159} & \textbf{0.7103} & 0.1333 & \textbf{0.0005} & \textbf{0.60} & \textbf{0.70} \\
		& PIQ & 0.0000 & \textbf{0.1250} & 1.00 & 7.8821 & 1.5097 & 0.9833 & 0.0099 & 0.00 & 0.00 \\
		\hline
		(T2, Case 1) & PM & 0.4735 & \textbf{0.1516} & 0.00 & 1.1669 & 5.3420 & 0.6967 & 0.0003 & 0.09 & 0.13 \\
		& AHuber & - & - & - & 2.4934 & 4.5658 & 0.8300 & 0.0333 & 0.00 & 0.01 \\
		& PWLAD & 0.1900 & 0.3350 & 0.00 & 1.5265 & 12.1056 & \textbf{0.0667} & 0.0736 & 0.00 & \textbf{0.80} \\ 
		& PWHL & \textbf{0.0550} & 0.3375 & \textbf{0.40} & \textbf{0.7215} & \textbf{0.7201} & 0.1333 & \textbf{0.0008} & \textbf{0.60} & 0.70 \\
		& PIQ & 0.2925 & 0.1981 & 0.00 & 1.3197 & 0.9598 & 0.4200 & 0.0057 & 0.00 & 0.21 \\
		\hline
		(T2, Case 2) & PM & 0.6960 & 0.2908 & 0.00 & \textbf{0.4641} & 13.7066 & 0.0400 & 0.0122 & 0.01 & \textbf{0.90} \\
		& AHuber & - & - & - & 1.3476 & 1.6083 & 0.1300 & 0.2593 & 0.00 & 0.86 \\ 
		& PWLAD & 0.0100 & 0.3488 & 0.80 & 1.3548 & 13.0757 & 0.1667 & 0.0605 & 0.00 & 0.70 \\ 
		& PWHL & 0.0050 & \textbf{0.1113} & 0.90 & 0.7640 & \textbf{0.7458} & 0.1000 & \textbf{0.0025} & \textbf{0.40} & 0.80 \\
		& PIQ & \textbf{0.0070} & 0.1268 & \textbf{1.00} & 1.8169 & 1.1137 & \textbf{0.4933} & 0.0062 & 0.00 & 0.15 \\
		\hline
		(T2, Case 3) & PM & 0.0010 & 0.1744 & 0.98 & 1.1692 & 5.0517 & 0.6967 & 0.0009 & 0.07 & 0.16 \\
		& AHuber & - & - & - & 1.8101 & 2.8766 & 0.8700 & 0.0291 & 0.00 & 0.13 \\  
		& PWLAD & 0.0000 & 0.3475 & 1.00 & 1.2355 & 12.5870 & \textbf{0.1667} & 0.0637 & 0.00 & \textbf{0.60} \\
		& PWHL & 0.0000 & \textbf{0.1225} & 1.00 & \textbf{0.8899} & \textbf{0.8877} & 0.2333 & \textbf{0.0015} & \textbf{0.40} & 0.50 \\ 
		& PIQ & 0.0000 & 0.1250 & 1.00 & 7.6467 & 1.4871 & 0.9800 & 0.0099 & 0.00 & 0.00 \\
		\hline
	\end{longtable}	
\end{footnotesize}	

\begin{footnotesize}
	\begin{longtable}{lcccccccccc}
		\caption{Comparison of PM\citep{kong18}, PWLAD, AHuber, PIQ and the proposed method PWHL
			under three cases in heterogeneous model for $c=0.3$.}\label{tab2-3}\\
		\hline
		(T, Case)  & Method & M & S & JD &  EE & EE$_{non}$  & FZR & FPR & SR & CR \\
		\hline
		\endfirsthead	
		\hline
		(T, Case)  & Method & M & S & JD &  EE & EE$_{non}$  & FZR & FPR & SR & CR \\
		\hline
		\endhead
		\hline
		\endfoot
		(T1, Case 1) & PM & 0.6440 & \textbf{0.0967} & 0.00 & 1.3131 & 2.8028 & 0.8900 & 0.0003 & 0.03 & 0.05 \\
		& AHuber & - & - &- & 2.2533 & 2.8404 & 0.9133 & 0.0272 & 0.00 & 0.02 \\ 
		& PWLAD & 0.1733 & 0.4157 & \textbf{0.20} & 1.8959 & 10.3069 & 0.3333 & 0.0758 & 0.00 & 0.00 \\  
		& PWHL & \textbf{0.1167} & 0.3371 & 0.00 & \textbf{0.8544} & \textbf{0.8515} & \textbf{0.2000} & \textbf{0.0008} & \textbf{0.50} & \textbf{0.50} \\
		& PIQ & 0.1667 & 0.2714 & 0.00 & 1.4224 & 1.1478 & 0.7000 & 0.0078 & 0.00 & 0.10 \\
		\hline
		(T1, Case 2) & PM & 0.7067 & 0.2929 & 0.00 & \textbf{0.4073} & 13.9575 & \textbf{0.0200} & 0.0117 & 0.01 & \textbf{0.94} \\ 
		& AHuber & - & - & - & 1.3664 & 1.4325 & 0.1567 & 0.2496 & 0.00 & 0.84 \\ 
		& PWLAD  & 0.0000 & 0.3771 & 1.00 & 1.1850 & 11.8871 & 0.0667 & 0.0476 & 0.00 & 0.80 \\  
		& PWHL & 0.0000 & 0.2743 & 1.00 & 0.5969 & \textbf{0.5385} & 0.0333 & \textbf{0.0048} & \textbf{0.40} & 0.90 \\
		& PIQ & 0.0000 & \textbf{0.2000} & 1.00 & 1.1496 & 0.9259 & 0.3667 & 0.0053 & 0.00 & 0.20 \\
		\hline
		(T1, Case 3) & PM & 0.4577 & 0.1014 & 0.53 & 5.2460 & 16.0726 & 0.6833 & 0.0284 & 0.03 & 0.07 \\ 
		& AHuber & - & - & -  & 1.7844 & 2.7623 & 0.8833 & 0.0260 & 0.00 & 0.11 \\  
		& PWLAD & 0.0000 & 0.3486 & 1.00 & 1.2354 & 11.6776 & 0.1333 & 0.0630 & 0.00 & 0.70 \\ 
		& PWHL & 0.0000 & \textbf{0.0771} & 1.00 & \textbf{0.7896} & \textbf{0.7524} & \textbf{0.1000} & \textbf{0.0040} & \textbf{0.20} & \textbf{0.80} \\
		& PIQ & 0.0000 & 0.2000 & 1.00 & 7.5050 & 1.4297 & 0.9967 & 0.0101 & 0.00 & 0.00 \\
		\hline
		(T2, Case 1) & PM & 0.7353 & \textbf{0.0983} & 0.00 & 1.3504 & 2.1895 & 0.9467 & 0.0001 & 0.00 & 0.01 \\
		& AHuber & - & - & - & 2.1923 & 2.5903 & 0.9233 & 0.0222 & 0.00 & 0.01 \\
		& PWLAD & 0.2000 & 0.2971 & 0.00 & 1.9560 & 9.9790 & 0.3667 & 0.0758 & 0.00 & 0.20 \\ 
		& PWHL & \textbf{0.1200} & 0.3229 & 0.00 & \textbf{0.8722} & \textbf{0.8608} & \textbf{0.1667} & \textbf{0.0030} & \textbf{0.30} & \textbf{0.50} \\
		& PIQ & 0.1783 & 0.2764 & 0.00 & 1.4734 & 1.1560 & 0.6767 & 0.0076 & 0.00 & 0.06 \\
		\hline
		(T2, Case 2) & PM & 0.6960 & 0.2923 & 0.00 & 0.4775 & 13.5605 & \textbf{0.0400} & 0.0119 & 0.00 & \textbf{0.89} \\
		& AHuber & - & - & - & 1.3675 & 1.4250 & 0.1800 & 0.2602 & 0.00 & 0.82 \\ 
		& PWLAD & 0.0167 & 0.3700 & 0.60 & 1.2984 & 9.3754 & 0.2667 & 0.0531 & 0.00 & 0.50 \\
		& PWHL & 0.0200 & 0.2729 & 0.50 & \textbf{0.8086} & \textbf{0.7593} & 0.1667 & 0.0093 & \textbf{0.40} & 0.60 \\
		& PIQ & \textbf{0.0023} & \textbf{0.2010} & \textbf{1.00} & 1.3824 & 1.0230 & 0.5100 & \textbf{0.0064} & 0.00 & 0.15 \\
		\hline
		(T2, Case 3) & PM & 0.4313 & 0.1180 & 0.55 & 5.1660 & 12.5500 & 0.7300 & 0.0263 & 0.01 & 0.05 \\
		& AHuber & - & - & - & 1.8339 & 2.8645 & 0.8500 & 0.0354 & 0.00 & 0.15 \\   
		& PWLAD & 0.0000 & 0.3000 & 1.00 & 1.6602 & 11.6240 & 0.3000 & 0.0615 & 0.00 & 0.20 \\
		& PWHL & 0.0000 & \textbf{0.0943} & 1.00 & \textbf{0.9771} & \textbf{0.9727} & \textbf{0.2333} & \textbf{0.0018} & \textbf{0.20} & \textbf{0.40} \\
		& PIQ & 0.0000 & 0.2000 & 1.00 & 7.2773 & 1.3854 & 0.9967 & 0.0101 & 0.00 & 0.00 \\
		\hline
	\end{longtable}	
\end{footnotesize}	

\newpage
\subsection{Breast cancer survival data analysis: Selected genes by PIQ}

\begin{figure}[h]
	\begin{center}
		\includegraphics[width=.7\textwidth]{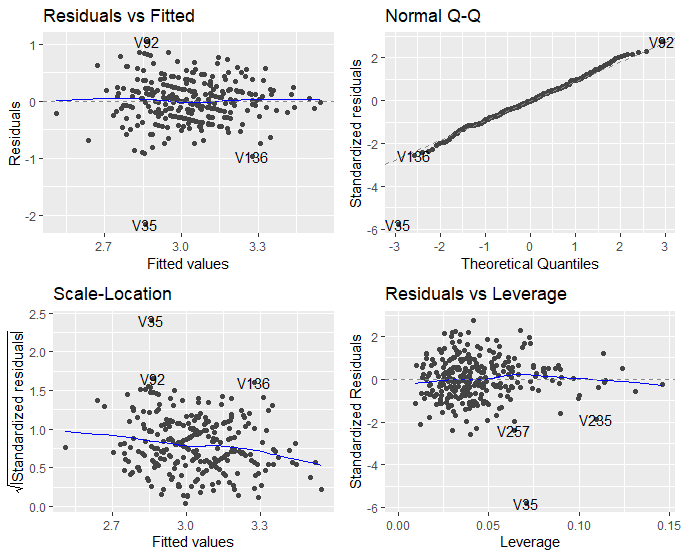}
	\end{center}
	\caption{Figures for the fitted PM model, including Residuals vs Fitted plot( upper left), QQ-plots of the residuals(upper right), Scale-Location plot(left below), and Residuals vs Leverage plot( right below).\label{fig2}}
\end{figure}

\begin{figure}[h]
	\begin{center}
		\includegraphics[width=.7\textwidth]{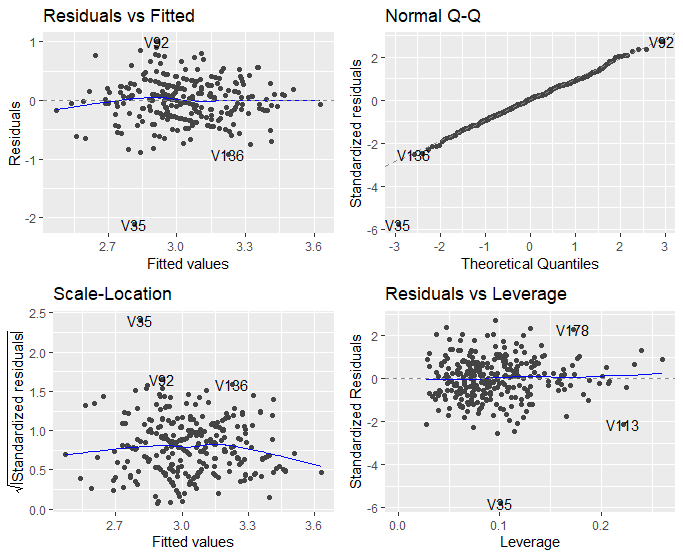}
	\end{center}
	\caption{Figures for the fitted PWHL model,  Residuals vs Fitted plot(upper left), QQ-plots of the residuals(upper right), Scale-Location plot(
		left below), and Residuals vs Leverage plot(right below).\label{fig3}}
\end{figure}

We present Scale-Location graphs in Figures \ref{fig2}--\ref{fig3} using PM\citep{kong18} and PHWL\citep{gao18} methods respectively in Section 4.2 to reveal endowed heterogeneity in breast cancer dataset.

\begin{footnotesize}
	\setlength{\tabcolsep}{1pt}
	\begin{longtable}{lccccc}
		\caption{Selected genes and corresponding coefficients by PIQ\citep{she22}}\\
		\hline
		\multicolumn{1}{c}{ Method}&
		\multicolumn{5}{c}{ Selected genes}\\
		\hline
		\endfirsthead	
		\hline
		\multicolumn{1}{c}{ Method}&
		\multicolumn{5}{c}{ Selected genes}\\
		\hline
		\endhead
		\hline
		\endfoot
		PIQ & NM.004701  & AF279865 & Contig50981.RC & NM.014364 & NM.019013 \\
		& 0.1193 & -0.1210  & -0.1607  & 0.3989  & 0.1162 \\
		& NM.006115 & Contig21406.RC & NM.004207 & NM.003258 & D14678\\
		& 0.1695  & -0.3908  & -0.0771  & -0.1637  & 0.0055 \\
		& Contig42342.RC & NM.004867 & NM.003022 & NM.014750 & NM.003295\\
		& -0.4071 & -0.0005  & 0.1178  & -0.1399 & -0.3401 \\
		& Contig46991.RC & NM.006633 & NM.003579 & NM.013999 & NM.004456\\
		& 0.0012  & -0.7045  & 0.1619  & 0.0435  & 0.2922\\ 
		& NM.001168 & M64497 & AK001166 & NM.001677 & Contig41413.RC\\
		& 0.1586  & -0.0435  & 0.2836  & -0.1665  & -0.0592 \\
		& NM.002989 & NM.001360 & Contig693.RC & NM.003740 & NM.014321\\
		& 0.0655  & -0.0471  & -0.1865  & 0.0737  & 0.1219 \\
		& NM.006829 & AL049949 & NM.000849 & Contig50855.RC & NM.006006\\
		& -0.2133  & 0.2657  & -0.0104  & 0.2568  & 0.2156 \\
		& NM.005321 & Contig54867.RC & NM.001333 & D25328 & NM.005181\\
		& 0.4464 & 0.0701   & 0.3512  & -0.0508  & -0.2180 \\
		& NM.006475 & NM.006101 & NM.001786 & AB033114 & AL133619\\
		& -0.2693  & -0.4053  & 0.0492  & -0.0124  & 0.0718 \\
		& AF070632 & NM.003108 & NM.004967 & NM.001809 & Contig778.RC\\
		& -0.1669  & 0.0498 & 0.1089  & -0.1313  & -0.1627\\
		& Contig45316.RC & NM.004217 & Contig38928.RC & NM.006386 & NM.018407\\
		& -0.2389  & -0.7691  & -0.1352  & -0.0388  & -0.1165 \\
		& Contig14284.RC & NM.000237 & NM.006096 & Contig1982.RC & NM.000561\\
		& 0.1441  & -0.0274  & 0.0662  & -0.0561  & 0.3435 \\
		& NM.005555 & NM.018087 & NM.017680 & Contig56390.RC & NM.002421\\
		& -0.3300  & 0.1749  & 0.2966  & -0.0514  & -0.2122 \\
		& NM.015417 & NM.001609 & NM.003613 & Contig55997.RC & Contig45537.RC\\
		& -0.2939  & -0.0601 & -0.1127   & 0.1561  & 0.0419 \\
		& NM.016441 & NM.016352 & NM.002193 & J03817 & Contig14658.RC\\
		& -0.2100  & 0.2450  & -0.2857  & -0.4401  & -0.5940 \\
		& NM.005410 & NM.004648 & AL117418 & NM.020974 & AL157502\\
		& 0.0142  & -0.2573 & 0.3066  & 0.0270  & 0.0893 \\
		& M96577 & Contig58107.RC & NM.002110 & AB007950 & Contig33741.RC\\
		& -0.3599  & 0.4888  & -0.1217  & -0.2309  & 0.1154 \\
		& AF052159 & Contig44191.RC & Contig38170.RC & Contig53307.RC & Contig52305.RC\\
		& 0.6305  & 0.0713  & -0.0020  & -0.1076  & -0.3699 \\
		& AF041429 & AL137707 & NM.000667 & X99142 & NM.004295\\
		& 0.5866  & -0.1885  & 0.0591  & 0.0711  & 0.2967 \\
		& Contig49279.RC & NM\.018265 & Contig57138.RC & NM.003956 & NM.000111\\
		& -0.1023  & -0.4290  & 0.1219  & -0.3726  & 0.0644 \\
		& NM.000320 & NM.002258 & NM.004780 & NM.021190 & NM.006341\\
		& -0.1890  & -0.3378  & 0.1757  & -0.2586  & 0.2687 \\
		& Contig7558.RC & NM.002773 & NM.005824 & NM.006623 & NM.000662\\
		& -0.1509  & 0.5994  & 0.4566  & 0.0991  & 0.1206 \\
		& NM.014292 & AF079529 & Contig50814.RC & AI147042.RC & NM.005505\\
		& 0.2004  & 0.0825  & -0.0656  & 0.1553  & 0.4700 \\
		& NM.006344 & Contig46597.RC & Contig24609.RC & NM.018455 & Contig57293\\
		& 0.0519  & -0.2397  & -0.0084  & -0.0138  & -0.8053 \\
		& K02403 & Contig20629.RC & Contig54993.RC & NM.001394 & NM.014669\\
		& 0.0206  & -0.0236  & 0.0286  & 0.1734  & 0.2616 \\
		& Contig55048.RC & Contig35897.RC & Contig41850.RC & NM.007088 & NM.001267\\
		& 0.8783  & -0.1354  & 0.6258  & 0.1555  & -0.0761 \\
		& NM.007106 & NM.004071 & X03084 & Contig55771.RC & NM.001395\\
		& 0.2924  & -0.0692  & 0.1073  & -0.0301  & -0.2094 \\
		& NM.003264 & Contig19284.RC & NM.003311 & NM.004315 & NM.004003\\
		& 0.0590  & 0.5288  & 0.0093  & 0.0332  & -0.0240 \\
		& NM.001844 & Contig41275.RC & Contig58301.RC & NM.000993 & NM.000269\\
		& -0.2070  & -0.2449  & 0.0126  & -0.1476  & 0.1216 \\
		& NM.002928 & NM.004095 & Contig43759.RC & NM.005101 & NM.005804\\
		& 0.0747  & -0.5664  & 0.2367  & -0.0303  & -0.0042 \\
		& Contig50153.RC & NM.002811 & NM.001605 & NM.004472 & Contig38493.RC\\
		& 0.3964  & 0.2208  & -0.2071  & 0.2331  & -0.0839 \\
		& NM.004265 & Contig50367 & Contig50950.RC & NM.000633 & Contig46.RC\\
		& -0.0360  & -0.3359  & 0.0049  & -0.2022  & 0.3200 \\
		& NM.000582 & AF035284 & Contig48756.RC & AF134404 & NM.003090\\
		& 0.0988  & -0.1608  & -0.0067  & 0.0512  & -0.3265 \\
		& NM.005689 & NM.014298 & NM.002805 & NM.014373 & AL080235\\
		& -0.0747 & -0.0117  & 0.2851  & -0.1198  & -0.3445 \\
		& Contig8818.RC & & & & \\
		& -0.2516 & & & & \\    
		\hline
		
	\end{longtable}	
\end{footnotesize}

\end{document}